\documentclass[11pt]{article}
\usepackage{amsmath, amssymb}
\usepackage[english]{babel}

\newtheorem{prop}{Proposition}[section]
\newtheorem{thm}[prop]{Theorem}
\newtheorem{defn}[prop]{Definition}
\newtheorem{lem}[prop]{Lemma}
\newtheorem{fact}[prop]{Fact}
\newtheorem{pr}[prop]{Property}
\newtheorem{cor}[prop]{Corollary}
\newtheorem{hyp}[prop]{Hypothesis}

\def\P{\Bbb P}

\def\E{\Bbb E}
\def\C{\Bbb C}
\def\e{\epsilon}
\newcommand{\RR}{\mathbb{R}}
\newcommand{\CC}{\mathbb{C}}
\newcommand{\NN}{\mathbb{N}}

\newcommand{\EE}{\mathbb{E}}

\newcommand{\PR}{\mathcal{P}(\RR)}

\newcommand{\MNC}{\mathcal{M}_N(\CC)}

\newcommand{\ppq}{\leqslant}
\newcommand{\pgq}{\geqslant}

\def\N{{\Bbb N}}
\def\prf{{\bf Proof.}}

\def\tr{{\mbox{tr}}}

\def\nn{\noindent}

\def\eps{\epsilon}

\def\diag{{\mbox{diag}}}

\def\ldp{{large deviation principle\ }}

\def\xx{\vrule height 0.7em depth 0.2em width 0.5 em}

\def\a{\alpha}
\def\b{\beta}
\def\d{\delta}
\def\e{\epsilon}
\def\g{\gamma}

\def\k{\kappa}
\def\l{\lambda}

\def\s{\sigma}
\def\t{\theta}

\def\z{\zeta}

\def\z{\zeta}

\def\D{\Delta}
\def\G{\Gamma}
\def\L{\Lambda}

\def\ra{\rightarrow}
\def\longra{\longrightarrow}

\def\Aa{{\mathcal A}}

\def\Ca{{\mathcal C}}
\def\Da{{\mathcal D}}

\def\Ia{{\mathcal I}}

\def\La{{\mathcal L}}

\def\Oa{{\mathcal O}}

\def\Ua{{\mathcal U}}

\def\mun{{\hat\mu^N}}

\def\lbc{\lbrace}

\def\lbk{\lbrack}
\def\rbk{\rbrack}
\def\nn{\noindent}
\def\part{\partial}

\def\ot{\otimes}

\def\ts{\times}
\def\Pa{{\mathcal{P}}}

\def\Ip{I^{\prime}}
\def\Ipp{I^{\prime \prime}}

\def\HE{H_{\mu_E}}
\def\KE{K_{\mu_E}}
\def\RE{R_{\mu_E}}

\def\QE{Q_{\mu_E}}
\def\IE{I_{\mu_E}}

\def\HM{H_{max}}
\def\Hm{H_{min}}
\def\lM{\l_{max}}
\def\lm{\l_{min}}

\def\ldp{large deviations principle}
\def\trace{\mathop{\rm tr}\nolimits}

\title{A Fourier view on the $R$-transform and related
asymptotics of spherical integrals}

\author{Alice Guionnet\thanks{Ecole Normale Sup\'erieure de Lyon,
Unit\'e de Math\'ematiques pures et appliqu\'ees,
UMR 5669,
46 All\'ee d'Italie, 
69364 Lyon Cedex 07, France. E-mail: aguionne@umpa.ens-lyon.fr.}, Myl\`ene Ma\"{\i}da\thanks{Universit\'e Paris X-Nanterre,
Laboratoire Modal-X, 200 av. de la R\'epublique, 92001 Nanterre Cedex, France. E-mail: mmaida@umpa.ens-lyon.fr}}

\date{}

\begin{document}

\renewcommand{\refname}{\Large{References}}

\maketitle

\begin{abstract}
We estimate the asymptotics of spherical
integrals of real symmetric or Hermitian matrices when the rank of one matrix
is much smaller than its dimension. We show that it is given in
terms of the $R$-transform of the spectral
measure of the full rank matrix
and give a new proof of the
fact that the $R$-transform is additive
under free convolution.  
These asymptotics also extend to the case where
one matrix has rank one but complex eigenvalue,
a result related with the analyticity 
of the corresponding spherical integrals. \\ %Ader

\end{abstract}

\section{Introduction}

\subsection{General framework and statement of the results}

In this article,  we consider the 
spherical integrals 
$$I_N^{(\beta)}(D_N,E_N):=
\int \exp\{N\trace(UD_NU^*E_N)\} dm_N^\beta(U),$$
where $m_N^{(\beta)}$ denote the Haar measure on the orthogonal
 group $\Oa_N$
when $\beta=1$ and  on the
 unitary group
${\mathcal U}_N$
when $\beta=2$, and $D_N$, $E_N$ are $N\ts N$
 matrices that we can assume diagonal without
loss of generality.
Such integrals are often called, in the physics literature,
Itzykson-Zuber or Harich-Chandra 
integrals. We do not consider the case $\b=4$ 
mostly to 
 lighten the
notations.\\
The interest for these objects goes back in particular to
the work of Harish-Chandra (\cite{HC1}, \cite{HC2}) 
who intended to define a notion of Fourier transform on Lie algebras.
They have been then extensively studied in the framework of so-called
matrix models that are related to the problem of enumerating maps
 (after \cite{ItZ}, it has been developed in physics
for example in \cite{ZZJ}, \cite{K2} or \cite{Ma},
in mathematics in \cite{Co} or \cite{G};
a very nice introduction to these links is provided in 
 \cite{ZV}). The asymptotics of the  
spherical integrals needed to solve matrix models were investigated
in \cite{GZ}. More precisely, when $D_N$, $E_N$
have $N$ distinct real eigenvalues $(\t_i(D_N),\l_i(E_N))_{1\ppq i\ppq N}$
 and the spectral measures
$\mun_{D_N}={1\over N}\sum\d_{\t_i(D_N)}$ and $\mun_{E_N}= {1\over N} \sum \d_{\l_i(E_N)}$
converge respectively to $\mu_D$ and  $\mu_E$, it is proved in Theorem 1.1 of \cite{GZ}
that
\begin{equation}\label{gzres}
\lim_{N\ra\infty}{1\over N^2}\log I_N^{(\beta)}(D_N,E_N)=I^{(\beta)}(\mu_D,\mu_E)
\end{equation}
exists under some technical
assumptions and a (complicated) formula for this limit is given.\\

In this paper, we investigate different asymptotics
of the spherical integrals, namely the case where one of the
matrix, say $D_N$, has  rank much smaller than
 $N$.\\

Such asymptotics were also already used in 
physics (see \cite{MPR}, where they consider
replicated 
spin glasses, the number of replica being there
the rank of $D_N$)
or stated for instance in \cite{Co}, section 1,
as a formal limit (the spherical 
integral being seen as a serie %Ader
 in $\theta$
when $D_N=\mbox{diag}(\theta,0,\cdots,0)$ 
whose coefficients are converging
as $N$ goes to infinity). However,
to our knowledge, there is no rigorous
derivation of this limit
available in the literature. We here
study this problem by use of
large deviations techniques. The proofs
are however rather different from 
those of \cite{GZ} ;
they rely on
large deviations for Gaussian
variables and not on their Brownian motion
interpretation and stochastic analysis as in \cite{GZ}.

\medskip

Before stating our results, 
we now introduce some notations and make a few remarks.\\
Let $D_N=\mbox{diag}(\theta,0,\cdots,0)$ have rank one
so that 
\begin{equation}\label{r1}
I_N^{(\b)}(D_N,E_N)=I_N^{(\b)}(\theta,E_N)=
\int e^{\theta N (UE_NU^*)_{11}} dm_N^{(\beta)}(U).
\end{equation}
\noindent Note that in general, in the case $\beta=1$, we will omit the superscript
$(\beta)$ in all these notations.\\

\noindent We make the following hypothesis :
\begin{hyp}
\label{hypsurEN}
\mbox{}
\begin{enumerate}
\item $\mun_{E_N}$ converges weakly towards a compactly 
supported measure $\mu_E$.
\item $\l_{\min}(E_N):=\min_{1\ppq i\ppq N}
\l_i(E_N)$ and $\l_{\max}(E_N):=
\max_{1\ppq i\ppq N} \l_i(E_N)$ converge respectively
to $\l_{\min}$ and $\l_{\max}$ which are finite.\\
\end{enumerate}
\end{hyp} 
Note that under Hypothesis \ref{hypsurEN},
the support of $\mu_E$, which we shall denote $\mbox{supp}(\mu_E)$,
 is included into  
$[\l_{\min},\l_{\max}]$.

Let us denote, 
for a probability measure $\mu_E$, its Hilbert 
transform by $\HE$   :
\begin{equation}
\label{Hnew}
\HE : \begin{array}[t]{rcl}
I_E:= \RR  \setminus \mbox{supp}(\mu_E) & \longra & \RR \\
z & \longmapsto & {\displaystyle \int \frac{1}{z-\l} d\mu_E(\l).}
\end{array}
\end{equation}
It is easily seen 
(c.f subsection \ref{ppn} for details)
that $\HE : I_E\ra \HE(I_E)$ 
is invertible, with inverse denoted $K_{\mu_E}$.
For $z\in \HE(I_E)$,
we set $R_{\mu_E}(z)=K_{\mu_E}(z)-z^{-1}$ 
to be the so-called $R$-transform of $\mu_E$.
In the case of the spectral measure $\hat \mu_{E_N}^N$
of $E_N$, we denote by $H_{E_N}$ its Hilbert transform 
given by $H_{E_N}(x) = \frac 1 N \tr(x - E_N)^{-1}=
\frac 1 N \sum_{i=1}^N(x - \l_i(E_N))^{-1} .$

The central result of this paper
can be stated as follows :
\begin{thm}
\label{sh}
Let $\b=1$ or $2$. If we assume that 
Hypothesis \ref{hypsurEN}.1 is satisfied %Ader
and that there is $\eps >0$ such that
\begin{equation}
\label{hypnorm}
\| E_N\|_\infty := \max\{|\l_{max}(E_N)|,|\l_{min}(E_N)|\} = O\left(N^{\frac 1 2 -\eps }\right),
\end{equation}
then for $\theta$ small enough
so that there exists $\eta>0$ so that
\begin{equation}\label{cond2}
\frac{2\theta}{\b}\in \bigcup_{N_0\ge 0}
\bigcap_{N\ge N_0}
H_{E_N}([\l_{\rm min}(E_N)-\eta,\l_{\rm max}(E_N)+\eta]^c),
\end{equation}
\begin{equation}\label{int1}
\IE^{(\b)}(\t):=\lim_{N\ra\infty}{1\over N}\log I_N^{(\b)}(\theta,
E_N)={\beta\over 2}\int_0^{2\theta\over \beta}
R_{\mu_E}(v) dv.
\end{equation}
Under Hypothesis \ref{hypsurEN}.2, \eqref{hypnorm} is obviously satisfied and
(\ref{cond2}) is equivalent to 
$$\frac{2\theta}{\b}\in H_{\mu_E}([\l_{min},\l_{max}]^c).$$
\end{thm}
This result is proved in section \ref{short} 
and appears in a way as a by-product of Lemma \ref{cont}. 
It raises several remarks
 and generalisations
that we shall investigate in this paper.\\
Note that in Theorems \ref{precise1}, \ref{rg1complexe}
and \ref{Fourier} hereafter
we consider the case $\b =1$, which requires simpler notations 
but every statement could be extended to the case $\b =2$.
The  main difference to extend these theorems to the case $\b=2$ is
that, following Fact \ref{gaussien},
it requires to deal with twice as much Gaussian 
variables, and hence to consider 
covariance matrices with twice bigger
dimension (the difficulty
lying then in showing that these matrices are positive definite).\\

The first question we can ask is how to precise the convergence \eqref{int1}. 
Indeed, in the full rank asymptotics, in particular
in the framework of \cite{GZ}, the second order term 
has not yet been rigorously derived. In our case, if
$d$ is  the Dudley distance  between measures 
(which is compatible
with the weak topology) given by
\begin{equation}\label{Dudley}
d(\mu, \nu) = \sup \left\{ \left|\int f d\mu - \int f d\nu\right| ; \, 
|f(x)|\textrm{ and } 
\left| \frac{f(x)-f(y)}{x-y}\right| \ppq 1, \forall x \neq y\right\},
\end{equation}
 we have

\begin{thm}
\label{precise1}
Assume Hypothesis \ref{hypsurEN}
and $$ d(\mun_{E_N},\mu_E)=o(\sqrt{N}^{-1}).$$
\indent Let $\t$ be such that 
${2\theta}\in H_{\mu_E}([\l_{min},\l_{max}]^c)$. 

\begin{itemize}
\item
If 
$\mu_E$ is not a Dirac measure at a single  point,
then, with $v=\RE(2\theta)$,
$$\lim_{N\ra\infty} e^{-{N}\left(\t v-{1\over 2N}\sum_{i=1}^N
\log(1-2\theta\l_i(E_N) +2\theta v)
\right)  }
I_N(\theta, E_N)= \frac{\sqrt{Z-4\t^2}}{\t \sqrt Z}, $$
with $\displaystyle Z:=  \int \frac{1}{(\KE(2\t)-\l)^2}d\mu_E(\l).$

\item
If $\mu_E=\d_e$
for some $e\in\RR$,
$$\lim_{N\ra\infty} e^{-{N}\theta e }
I_N(\theta, E_N)= 1. $$
\end{itemize}
\end{thm}

This theorem gives the second order term for the convergence
given in Theorem \ref{sh} above. Indeed,
with ${2\theta}\in H_{\mu_E}([\l_{min},\l_{max}]^c)$,
under Hypothesis \ref{hypsurEN}.2, there exists  (c.f.
(\ref{boundv}) for details)
$\eta(\theta)>0$ so that for $N$ large enough
$$1-2\theta\l_i(E_N) +2\theta v>\eta(\theta).$$
Therefore, there exists a finite constant $C(\theta)\le
(\eta(\theta)^{-1}+|\log(\eta(\theta))|)$ such that
for $N$ sufficiently large
\begin{multline*}
\left|{1\over 2N}\sum_{i=1}^N
\log(1-2\theta\l_i(E_N) +2\theta v)-
{1\over 2}\int \log(1-2\theta \l+2\theta v)d\mu_E(\l)\right|\\
\le C(\theta)d\left({1\over N}\sum_{i=1}^N \d_{\l_i(E_N)},
\mu_E\right),
\end{multline*}
where $d$ is  the Dudley distance.  \\
Moreover, with $v=\RE(2\theta)$, it is easy to see that
$$\theta v-{1\over 2}\int \log(1-2\theta \l+2\theta v)d\mu_E(\l)
={1\over 2}\int_0^{2\theta} R_{\mu_E}(u)du,$$
showing how  Theorem \ref{precise1} 
relates with Theorem \ref{sh}.\\

%Ader :ce n'est pas ce que je voulais, voir plus loin.

\medskip

Another remark is that Theorem \ref{sh} can be seen as giving 
an interpretation
of the primitive of the  $R$-transform
$\RE$ as a Laplace transform 
of $(UE_NU^*)_{11}$ for large $N$ and for compactly supported
probability measures $\mu_E$.\\
\indent  A natural question is to wonder
whether it can be extended to
the case where $\theta$ is complex,
to get an analogy with the Fourier
transform that seems to have originally motivated Harish-Chandra.
In the case
of the different asymptotics studied in \cite{GZ},
this question is
open : in physics, formal analytic extensions
of the formula obtained for Hermitian
matrices to any matrices
are commonly used, but S. Zelditch \cite{Zel}
found that such an extension
could be false by exhibiting
counter-examples. In the context of
the asymptotics we consider here,
we shall however see that
this extension is valid for $|\theta|$ small 
enough. Note that, as far as $\mu_E$ is compactly supported,
$R_{\mu_E}$ can be extended analytically 
at least in a complex neighborhood of the origin
(see Proposition \ref{Hcomp} for further details).

\begin{thm}
\label{rg1complexe}
Take $\b=1$ and 
assume that $({E_N})_{N\in\NN}
$ is a uniformly bounded sequence
of matrices satisfying  Hypothesis \ref{hypsurEN}.1
where $\mu_E$ is not  a Dirac mass.\\
Assume furthermore that
$d(\mun_{E_N},\mu_E)=o(\sqrt{N}^{-1})$, where $d$ is 
the Dudley distance defined by \eqref{Dudley}.\\
Then, there exists an $r >0$ such that,
 for any $\theta\in\C$, such that $|\theta| \ppq r$,
$$\lim_{N\ra\infty}{1\over N}\log 
I_N(\theta, E_N)=\theta v(\theta)
-{1\over 2}\int\log(1 +2\theta
v (\theta)- 2 \t \l)d\mu_E(\l),
$$
where $\log(.)$ is the main branch of the logarithm in $\CC$ 
and $v(\t) = \RE(2\t)$. More precisely, we prove that
for $\theta$ in a small complex neighborhood  of 
the origin, 
$$\lim_{N\ra\infty} e^{-{N}\left(\t v-{1\over 2N}\sum_{i=1}^N
\log(1-2\theta\l_i(E_N) +2\theta v)
\right)  }
I_N(\theta, E_N)= \frac{\sqrt{Z-4\t^2}}{\t \sqrt Z}, $$
with $\displaystyle Z:=  \int \frac{1}{(\KE(2\t)-\l)^2}d\mu_E(\l).$

\end{thm} %Ader
It is not hard to see that the above convergence
is uniform in a small complex neighborhood  of 
the origin. Consequently, there exists $\theta_0>0$, 
$N_0\in\N$,
such that for $|\theta|\le \theta_0$,
for all $N\ge N_0$, 
$f_N(\t) := \frac 1 N \log I_N(\t, E_N)$, 
is bounded from above and below.
Moreover, under Hypothesis \ref{hypsurEN}, the $f_N$'s are 
holomorphic and 
uniformly bounded. 
Therefore, by Cauchy's formula
$$\partial^{(n)} f_N|_{z=0}=-
{1\over 2\pi i} \int_{|z|=\theta_0/2}
\frac{f_N(z)}{z^{n+1}}dz$$
insures with dominated convergence theorem's that
for all $n\in\N^*$,
$$\lim_{N\ra\infty} \partial^{(n)} f_N|_{z=0}=\partial^{(n)} 
f|_{z=0}= 2^{k-1} \partial^{(n-1)} \RE|_{z=0}$$
with
$f(\theta)=
\theta v(\theta)
-{1\over 2}\int\log(1 +2\theta
v (\theta)- 2 \t \l)d\mu_E(\l))$.
Hence, we give a new proof 
of  B.~Collins'
result \cite{Co} (here in the orthogonal setting
rather than in the unitary one)
and validate the strategy, commonly used in physics,
of computing $f$ to calculate $\lim_{N\ra\infty} \partial^{(n)}
 f_N|_{z=0}$.

Note that the case $\mu_E=\d_e$ 
is trivial if we assume additionnally Hypothesis \ref{hypsurEN}.2 with $\l_{\min}$ and $\l_{\max}$
the edges of the support
of $\mu_E$ since then $\max_{1\le i\le N}|\l_i-e|$
goes to zero with $N$ which entails

$$\lim_{N\ra\infty}{1\over N}\log 
I_N(\theta, E_N)=\theta e$$
for all $\t$ in $\CC$.

The proof of Theorem \ref{rg1complexe} will be more involved than the
real case treated in sections \ref{short} and \ref{tcl}
and the difficulty lies of course in
the fact that the integral is now oscillatory,
forcing us to control
more precisely the deviations
in order to make sure that
the term of order one in the 
large $N$ expansion does not vanish.
This is the object of section \ref{C}.\\

Once the view of spherical integrals
as Fourier transforms has been justified by the extension
 to the complex plane, a second natural question
is to wonder whether we can use it to see
that
the $R$-transform is additive
under free convolution.
Let us make some reminder about free probability :
in this set up,  the notion of freeness 
replaces the standard
notion of independence and 
the R-transform is analogous to the logarithm of the Fourier 
transform of a measure. Now, 
it is well known that
the log-Laplace (or Fourier) transform 
is additive under convolution i.e. for any probability
measures $\mu$, $\nu$ on $\RR$ (say compactly
supported to simplify), any $\l\in\RR$, (or $\C$)
$$\log \int e^{\l x}d\nu*\mu(x)=
\log\int e^{\l x} d\mu(x)+\log\int e^{\l x} d\nu(x).$$
Moreover, this property, if it holds
 for $\l$'s in a neighbourhood
of the origin,
 characterizes uniquely 
the convolution.
Similarly, if we denote $\mu\boxplus\nu$
the free convolution of two compactly supported probability
measures on $\RR$, it is uniquely
described by the fact that
$$R_{\mu\boxplus \nu}(\l)=R_\mu(\l)+R_\nu(\l)$$
for sufficiently small $\l$'s.
Theorem \ref{sh} provides an interpretation of this result.
Indeed, Voiculescu \cite{Vo1}
proved that if $A_N,B_N$
are two  diagonal matrices
with 
spectral measures converging towards $\mu_A$
and $\mu_B$ respectively, with uniformly bounded
spectral radius,  then
the spectral measure of $A_N+UB_N U^*$
converges, if $U$ follows $m_N^{(2)}$,
towards $\mu_A\boxplus\mu_B$. This result
extends naturally to the case where $U$ follows $m_N^{(1)}$
(see \cite{BP} Theorem 5.2 for instance).
Therefore, it is natural to expect
 the following result :
\begin{thm}
\label{Fourier}
Let $\beta=1$,  $(A_N,B_N)_{N\in\N}$ be a sequence of
 uniformly bounded real diagonal matrices
and 
 $V_N$ following $m_N^{(1)}$.

\begin{enumerate}
\item Then
\begin{equation}
\label{Conc}
\lim_{N\ra\infty}\left( \frac 1 N \log I_N(\theta, A_N+V_NB_NV_N^*) 
-\int  \frac 1 N \log I_N(\theta, A_N+V_NB_NV_N^*) d
m_N^{(1)}(V_N)\right)=0\mbox{ a.s.}
\end{equation}
\item If additionnally the spectral 
measures of $A_N$ and $B_N$ converge respectively to
$\mu_A$ and $\mu_B$ fast enough (i.e. such that
$d(\hat\mu_{A_N}, \mu_A) +d(\hat\mu_{B_N}, \mu_B) = o(\sqrt N^{-1})$)
and $\mu_A$ and $\mu_B$ are not Dirac masses at a point, then, 
for any $\theta$ small enough,
\begin{equation}
\label{Conc2} 
\lim_{N\ra\infty}   \frac 1 N \log I_N(\theta, A_N+V_NB_NV_N^*) 
=\lim_{N\ra\infty}   \frac 1 N \log I_N(\theta, A_N)+
\lim_{N\ra\infty}  \frac 1 N \log I_N(\theta, B_N) \mbox{ a.s} .
\end{equation}
\end{enumerate} 
\end{thm}

Then the additivity of the 
$R$-transform (cf. Corollary \ref{additive}) is a 
direct consequence of 
this result together with the continuity of the spherical 
integrals with respect to the empirical measure 
of the full rank matrix (which will be shown in Lemma \ref{cont}).\\
Note that the case where $\mu_A$ or $\mu_B$ are
Dirac masses is trivial if we assume that
the edges of the spectrum of $A_N$ or
$B_N$ converge towards this point.
The general case could be handled as
well but, since it has no motivation 
for the $R$-transform (for which we
can always assume that the above condition
holds, see Corollary \ref{additive}),
we shall not detail it.
Section \ref{add} will be devoted to the proof
of this theorem  which decomposes 
mainly in two steps : to get the first point,
we establish a result of concentration
under $m_N^{(1)}$
 that will give us \eqref{Conc}; then to prove the second point
 once we have the first one it is enough
to consider the expectation of $ \frac 1 N  \, \log I_N
 (\theta, A_N+V_NB_NV_N^*)$
and
if one assumes
that
\begin{multline}\label{interch}
\lim_{N \ra \infty} {1\over N}\int\left(\log \int e^{\theta
N(UA_NU^*+UV_NB_NV_N^* U^*)_{11}}
dm_N^{(1)}
(U)\right)d m_N^{(1)}
(V)\\
= \lim_{N \ra \infty}
{1\over N}\log\int  \int 
e^{\theta N(UA_NU^*+UV_NB_NV_N^* U^*)_{11}}
dm_N^{(1)}
(U)d m_N^{(1)}
(V)
\end{multline}
 the equality (\ref{Conc2})  follows 
from the observation that
the right hand side equals $N^{-1}  \, \log I_N
(\theta, A_N)+
 N^{-1}  \, \log I_N
(\theta, B_N)$. \\
\indent 
Note that equation (\ref{interch}) is rather typical to what should be
expected for  disordered particles systems in the
high temperature regime and indeed  our proof 
 follows some very smart ideas 
of Talagrand that he developed in the context
of Sherrington-Kirkpatrick
model of spin glasses at high temperature (see \cite{Tal}).
This proof is however rather technical
because the required control on the $L^2$ norm 
of the partition function  
is based on the study of second order corrections
of replicated systems which generalizes Theorem \ref{precise1}.\\

The next question, that we will
actually 
tackle in section \ref{rang1}, 
deals with the understanding of
the limit (\ref{int1})
for all the values of $\theta$.
We  find
the following result
\begin{thm}
\label{expreel} Let $\beta=1$ or $2$.
Assume $\mun_{E_N}$ satisfy Hypothesis \ref{hypsurEN}.\\
If we let $\displaystyle \Hm := \lim_{z \uparrow \lm} \HE(z)$
and  $\displaystyle \HM := \lim_{z \downarrow \lM} \HE(z)$,
then 
$$\lim_{N\ra\infty}{1\over N}\log I_N^{(\b)}(\theta,E_N)=
\IE^{(\b)}(\t) = 
\theta v(\theta) -{\beta\over 2}\int \log \!\Bigl(1+{2\over\beta}
\theta v(\theta) -{2\over\beta}\theta \l\Bigr)\!d\mu_E(\l)
$$
with
\begin{eqnarray*}
v(\theta)
=\left\{ \begin{array}[c]{ll}
R_{\mu_E}\bigl({2\over\beta}
\theta\bigr) & \textrm{ if } \Hm \ppq \frac{2\t}{\b}  \ppq \HM \\
\lM-{\beta\over 2\theta} 
 &  \textrm{ if }  \frac{2\t}{\b}  > \HM \\
\l_{\rm min}-{\beta\over 2\theta}
, &  \textrm{ if }  \frac{2\t}{\b}  < \Hm. 
\end{array}
\right.
\end{eqnarray*}
\end{thm}
Note here that the values of $\l_{\min}$ and $\l_{max}$
do affect the value of the limit of spherical integrals
in the asymptotics we consider here,
contrarily to what happens in the full rank asymptotics 
considered in \cite{GZ}.

As a consequence of Theorem \ref{expreel}, we can see that
there are two phase transitions at $\HM \b/2$ and $\Hm \b/2$
which are of second order in general
(the second derivatives of $\IE(\t)$
being discontinuous at these
points, except when $\l_{\max}\HE^\prime(\l_{\max})=1$ (or similar
equation with $\l_{\min}$
instead of $\l_{\max}$),
in which case the 
transition is of order 3). These
transitions can
in fact be  characterized by the asymptotic
behaviour of $(UE_NU^*)_{11}$ under
the Gibbs measure
$$d\mu_N^{\beta,\theta}(U)={1\over I_N^{(\beta)}(\theta,E_N)}
e^{N\theta (UE_NU^*)_{11}}dm_N^{(\beta)}(U).$$
For $\t \in \left[{\Hm \b \over 2},{\HM \b \over 2}\right]^c$,
$(UE_NU^*)_{11}$  saturates 
and converges $\mu_N^{\beta,\theta}$-almost surely 
towards $ \l_{max}-{\beta\over 2\theta}$
(resp. $\l_{\min}-{\beta\over 2\theta}$).
Hence, up to 
a small component of 
norm of order $\theta^{-1}$,
with high probability, the first column 
vector $U_1$ of $U$ will
align on the eigenvector
corresponding
to either the
smallest or the largest 
eigenvalue of $E_N$, whereas for smaller
$\t$'s, $U_1$ will
prefer to charge all the
eigenspaces of $E_N$.\\

Another natural question is to wonder what happens
when $D_N$ has not rank one but  rank
 negligible compared to $N$.
It is not very hard to see that in the case where
all the eigenvalues of $D_N$ are small enough
(namely when they all lie inside $\HE([\lm, \lM]^c)$),
we find that the spherical integral approximately factorizes
into a product of integrals of rank one.
More precisely, 

\begin{thm}
\label{rgfini}
Let $\beta=1$ or $2$.
Let $D_N = \mbox{diag}(\t_1^N, \ldots, \t_{M(N)}^N, 0, \ldots, 0)$
with $M(N)$ which is $o(N^{\frac 1 2 -\e})$ for some $\e >0$.
 Assume that $\mun_{E_N}$
fulfills Hypothesis \ref{hypsurEN}.1,
that $||E_N||_\infty=o({N}^{\frac 1 2 -\varepsilon} )$ 
for some $\varepsilon>0$ and that 
there exists $N_0\in\N$ and
 $\eta >0 $ such that, for all $N\ge N_0$ and $i$
from $1$ to $M(N)$,  
$\frac{2\t_i^N}{\b} \in H_{E_N}([\lm(E_N)- \eta, \lM(E_N) +\eta]^c)
$. \\
 Then, if 
$\displaystyle \frac {1}{ M(N)} \sum_{i=1}^{M(N)} \d_{\t_i^N}$
converges weakly to $\mu_D$, 
$$ \IE^{(\beta)}(D) 
:= \lim_{N \ra \infty} \frac{1}{N M(N)} \log I_N^{(\beta)}(D_N, E_N)$$
exists and is given by 
\begin{equation}
\label{sip}
\IE^{(\beta)}(D)  = \lim_{N \ra \infty} \frac{1}{M(N)} \sum_{i=1}^{M(N)} \IE^{(\beta)}(\t_i^N) = \int  \IE^{(\beta)}(\t) d\mu_D(\t) .
\end{equation}
\end{thm}

This will be shown at the end of section \ref{short}, the
proof being very similar to the case of rank one.
It relies mainly on Fact \ref{gaussien} hereafter
and comes from the fact that in such asymptotics the $M(N)$
first column vectors of an  orthogonal 
or unitary matrix distributed according
to the Haar measure behave approximately 
like independent vectors uniformly distributed 
on the sphere. This can be compared with
the very old result 
of E.~Borel \cite{Bor} which says that one entry
of an orthogonal matrix distributed according
to the Haar measure behaves like a Gaussian variable.
That kind of considerations finds continuation for example in a recent work of
A.~D'Aristotile, P.~Diaconis and C.~M.~Newman
\cite{ADN} where they consider a number of
element of the orthogonal group going to infinity
not too fast with $N$.
In the same direction, one can also mention the recent work
of T.~Jiang \cite{Ji} where he shows that the entries of the first
$O(N/\log N)$ columns of an Haar distributed 
unitary matrix can be simultaneously approximated 
by independent standard normal variables.\\
\indent Recently,  P.~\'Sniady could prove %Ader
by different techniques that the asymptotics 
we are talking about extend to 
$M(N) = o(N)$.\\
\indent Of course we would like to generalize also the full asymptotics
we've got in Theorem \ref{expreel} to the set up of finite rank
i.e. in particular consider the case where some (a $o(N)$ number)
of the eigenvalues of $E_N$ could converge away from the support.
It seems to involve not only the deviations of $\lM$ but those
of the first $M$ ones when the rank is $M$.
As it becomes rather
 complicate and as the proof is already rather involved in rank one, 
we postpone this issue to further research.

\medskip

To finish this introduction,
we also want to mention that
the results we've just presented give  (maybe) less obvious
relations between the $R$-transform 
and Schur functions or vicious
walkers.  \\
Indeed, if $s_\l$ denotes the Schur function associated
with a Young tableau $\l$ (cf. \cite{Sag} for more details),
then, it can be checked (cf. \cite{GM1} for instance)
that
$$
s_{\lambda}(M) = I_N^{(2)} \left(\log M, 
\frac l N  \right) 
\Delta \left({l\over N}\right)
 \frac{\Delta(\log M)}{\Delta(M)}
$$
with $l_i=\l_i+N-i$, $1\le i\le N$ and $\D(M)=\prod_{i<j}(M_i-M_j)$
when $M=\diag(M_1,\cdots,M_N)$.
Thus, our results also give the asymptotics
of Schur functions when $N^{-1}\d_{N^{-1}(\l_i+N-i)}$ converges towards 
some compactly supported probability measure $\mu$. For instance, 
 Theorem \ref{sh} implies that for $\t$ small enough
$$
\lim_{N\ra\infty}{1\over N}\log \left( {\prod_{i>j}(N^{-1} ({\l_j-j-\l_i+i}))^{-1}} s_{\lambda}(e^{\theta},1,\dotsc, \dotsc,1)
\right)
=\int_0^\theta R_\mu(u)du +\log ({\theta}({ e^\theta-1})^{-1}).$$

Such asymptotics should
 be  more directly related with  the
combinatorics of the symmetric group
and more precisely with non-crossing partitions
which play a key role in free convolution.\\
\indent On the other hand,
it is also known that spherical integrals are related
with the density kernel of vicious walkers,
that is Brownian motions conditionned to
avoid each others, either by using the
fact that the eigenvalues of the Hermitian Brownian motion 
are described by such vicious walkers (more commonly named
in this context Dyson's Brownian motions) or by applying
directly the result of Karlin-McGregor \cite{karlin}. 
Hence, the  study of the asymptotics of spherical integrals 
we are considering allows 
to estimate this density kernel
when $N-1$ vicious walkers start at the origin,
the last one starting at $\theta$
and at time  one reach $(x_1,\dotsc,x_N)$
whose empirical distribution approximates
a given compactly supported probability measure.

\subsection{Preliminary properties and notations}
\label{ppn}
Before going into the proofs themselves, we gather here 
some material and notations that will be useful throughout the paper.

\subsubsection{Gaussian representation of Haar measure}

\indent In the different cases we will develop, the first step 
will be always the same : we will represent the column vectors 
of unitary or orthogonal matrices distributed according 
to Haar measure via Gaussian vectors. To be more precise, 
we recall the following fact :\\

\begin{fact}
\label{gaussien}
Let $k \ppq N$ be fixed.\\
$\bullet$ Orthogonal case.\\
Let $U=(u_{ij})_{1 \ppq i,j \ppq N}$ be a random orthogonal matrix distributed according to $m_N^{(1)}$, the Haar measure on $\Oa_N$. 
Denote by $(u^{(i)})_{1\ppq i \ppq N}$ the column vectors of~$U$.\\
Let $(g^{(1)}, \ldots, g^{(k)})$ be $k$ independent standard Gaussian vectors in $\RR^N$ 
and let $(\tilde g^{(1)}, \ldots , \tilde g^{(k)})$ the vectors obtained from  $(g^{(1)}, \ldots, g^{(k)})$
by the standard Schmidt orthogonalisation procedure.\\
Then it is well known that
$$  (u^{(1)}, \ldots, u^{(k)}) \sim \left( \frac{\tilde g^{(1)}}{\left\| \tilde g^{(1)}\right\|}, \ldots,
\frac{\tilde g^{(k)}}{\left\| \tilde g^{(k)}\right\|} \right),$$
where $\|.\|$ denotes the Euclidean norm in $\RR^N$ and the equality $\sim$ means that the two $k \times N$-matrices
have the same law.\\
$\bullet$ Unitary case.\\
With the same notations, let $U$ be distributed according to
$m_N^{(2)}
$,  the Haar measure on $\Ua_N$. 
Let $(g^{(1), R}, \ldots, g^{(k),R}, g^{(1), I}, \ldots, g^{(k),I})$ be $2k$ independent standard Gaussian vectors in $\RR^N$
and let \linebreak
$(\tilde G^{(1)}, \ldots, \tilde G^{(k)})$ be the $k$ vectors obtained from $(g^{(1),R}+ ig^{(1), I}, \ldots, g^{(k),R}+ ig^{(k), I})$
by the standard Schmidt orthogonalisation procedure with respect to the usual scalar product in $\CC^N$.\\
Then we get that 
$$  (u^{(1)}, \ldots, u^{(k)}) \sim \left( \frac{\tilde G^{(1)}}{\left\| \tilde G^{(1)}\right\|}, \ldots,
\frac{\tilde G^{(k)}}{\left\| \tilde G^{(k)}\right\|} \right),$$
where  $\|.\|$ denotes the usual norm in $\CC^N$.
\end{fact} 

Note that heuristically, the above representation in terms of Gaussian vectors allows us to understand why the limit in the finite rank case behaves as
a sum of functions of each of the eigenvalues of $D_N$. Indeed, in high dimension, we know that a bunch of $k$ 
(independent of the dimension) Gaussian vectors are almost orthogonal one from another so that the orthogonalisation 
procedure let them almost independent. \\

\subsubsection{Some properties of the Hilbert and the R-transforms
of a compactly supported probability measure on $\RR$}

 Let $\l_{\min}(E)$ and $\l_{\max}(E)$ be
the edges of the support of $\mu_E$.
For all $\l_{\min} \ppq \l_{\min}(E)$ and $\l_{\max} \pgq \l_{\max}(E)$,
let us  denote by $\displaystyle \Hm := \lim_{z \uparrow \lm} \HE(z)$
and  $\displaystyle \HM := \lim_{z \downarrow \lM} \HE(z)$,
where $\HE$ was defined in \eqref{Hnew}.\\

We sum up the properties of $\HE$ that will be useful for us 
in the following
\begin{pr}
\label{pteH}
\mbox{:}
\begin{enumerate}
\item $\HE$ is decreasing and positive on $\{z > \lM\}$
 and decreasing and  negative on $\{z < \lm \}$.
\item Therefore, $\Hm$ exists in $\RR_-^* \cup \{ - \infty\}$
and $\HM$ exists in $\RR_+^* \cup \{ + \infty\}$.
\item $\HE$ is bijective from $I=\RR\backslash [\l_{\rm min},\l_{\rm
max}]$ onto its image 
$I^\prime := ]\Hm, \HM[ \setminus \{0\}$.
\item $\HE$ is analytic on $I$ and its derivative never cancels on $I$.
\end{enumerate}
\end{pr}
The third point of the property above allows the following
\begin{defn}
\label{KRreel}
\mbox{:}
\begin{enumerate}
\item $\KE$ is defined on $\Ip$ as the functional inverse of $\HE$.
\item  $\Ip$ does not contain $0$ so that, on $\Ip$, we can define
$\RE$ given by $\RE(\g) = \KE(\g) - \frac 1 \g $ for any $\g \in \Ip$.
\end{enumerate}  
\end{defn}
We will need to consider the inverse $\QE$ of $\RE$.
To define it properly, we have to look more carefully at the properties 
of  $\RE$. We have :
\begin{pr}
\label{Rreel}
\mbox{:}
\begin{enumerate}
\item $\KE$ and  $\RE$ are analytic (and in particular continuously differentiable) on $\Ip$.
\item  $\RE$ is increasing and its derivative never cancels.
\item $ \displaystyle \lim_{\g \ra 0^-}  \RE(\g) =   \lim_{\g \ra 0^+}  \RE(\g) = m := \int \l d\mu_E(\l).$
\item  $\RE$ is bijective from $\Ip$ onto its image $\Ipp := \displaystyle 
\left]  \lm - \frac 1 \Hm, \lM -  \frac 1 \HM\right[ \setminus \{m\}$
so that we can define its inverse $\QE$ from $\Ipp$ to $\Ip$. Moreover, $\QE$ is differentiable
on $I^{\prime \prime}$.
\end{enumerate} 
\end{pr}
The proof of these properties is easy and left to the reader.\\

\noindent The following property deals with the behaviour of these functions
on the complex plane. A proof of it can be found for example in 
\cite{VoiR}. We first extend the definition of the Hilbert transform, 
that we denote again $\HE$ by
\begin{equation}
\HE : \begin{array}[t]{rcl}
 \CC  \setminus \mbox{supp}(\mu_E) & \longra & \CC \\
z & \longmapsto & {\displaystyle \int \frac{1}{z-\l} d\mu_E(\l).}
\end{array}
\end{equation}
\mbox{} \\
\begin{pr}
\mbox{:}
\label{Rcomp}
\begin{enumerate}
\item There exists a neighbourhood $\Aa$
of $\infty$ such that $\HE$ is bijective
from $\Aa$ into $\HE(\Aa)$, which is a 
neighbourhood of $0$.
\item We denote by $\KE^{(c)}$ its functional inverse on $\HE(\Aa)$
and $\RE^{(c)}$ is given by $\RE^{(c)}(\g) = \KE^{(c)}(\g) - \frac 1 \g $
 for any $\g \in \HE(\Aa)$ (that does not contain $0$).
\item $\RE^{(c)}$ is analytic and coincides with $\RE$
on $\Ip \cap \HE(\Aa)$. Therefore, we denote it again $\RE$.
\end{enumerate} 
\end{pr}

\medskip

Note that throughout the paper, we will denote $\l_i := \l_i(E_N)$,
 $\t_i := \t_i(D_N)$ (and even $\t$ will denote $\t_1(D_N)$ in
the case of rank one) and   denote in short  
$H_{E_N}(x) = \frac 1 N \tr (x-E_N)^{-1}.$\\ %Ader

We now state the following property, which will
be useful in the proof of Theorem \ref{rg1complexe} :
\begin{prop}
\label{Hcomp}
If $(E_N)_{N \in \NN}$ is uniformly bounded and satisfying
Hypothesis \ref{hypsurEN}.1, there exists $r >0$
such that, for any $\t \in \CC$ such that $|\t| \ppq r$,
there is a solution of 
$$ H_{E_N}\left(\frac{1}{2\t} + v_N(\t)\right)= 2\t,$$
such that $v_N(\t) \xrightarrow[N \ra \infty]{} \RE(2\t).$
\end{prop}

\textbf{Proof of Proposition \ref{Hcomp} :}
Let $\Aa_N$ be a neighbourhood of $\infty$
on which $H_{E_N}$ is invertible
($\Aa_N$ can be given as $\{z/ |z| > R_N\},$ for some $R_N$).
For any $\eta >0$, we denote by
$\Aa_N^\eta := \{ x \in \Aa_N / d(x, \Aa_N^c) \pgq \eta\}.$
Let $\t$ be such that there exists $\eta >0$ such that
$2\t \in \bigcup_{N_0 \pgq 0}  \bigcap_{N \pgq N_0} H_{E_N}(\Aa_N^\eta),$
we take $v_N(\t)$ the unique solution in $\Aa_N^\eta - (2\t)^{-1}$
of 
$$ H_{E_N}\left(\frac{1}{2\t} + v_N(\t)\right)= 2\t.$$
Since, for all $\l \in \bigcup_{N_0 \pgq 0}  \bigcap_{N \pgq N_0} 
\mbox{supp}(\mun_{E_N})$, the application $z \mapsto (z-\l)^{-1}$
is continuous bounded on 
 $\bigcup_{N_0 \pgq 0}  \bigcap_{N \pgq N_0} \Aa_N^\eta$,
under Hypothesis \ref{hypsurEN}.1, $v_N(\t)$ converges to $\RE(2\t).$\\
Furthermore, the fact that $(E_N)_{N \in \NN}$ is uniformly bounded
ensures that we can choose the $\Aa_N$'s such that there exists $r >0$
such that 
$\bigcup_{N_0 \pgq 0}  \bigcap_{N \pgq N_0} H_{E_N}(\Aa_N^\eta)
\supset \{\t / |\t| \ppq r\}.$ \hfill\xx

%%%%%%%%%%%%%%%%%%%%%%%%%%%%%%%%%%%%%%%%%%%%%%%%%%%%%%%%%%%%%%%%%%%%

\section{Proof of Theorems \ref{sh}, \ref{rgfini}  and related results
}
\label{short}
Before going into more details, let us state and prove a lemma 
which deals with the continuity 
of $I_N$ and its limit. We state here
a trivial continuity in the finite rank
matrix but also a weaker continuity
result in the spectral measure
of the diverging rank
matrix, on which the proof of Theorem \ref{sh}
is based.

\begin{lem}
\label{cont} 
\begin{enumerate}
\item  For any $N\in\N$, any sequence of matrices $(E_N)_{N\in\N}$
with spectral radius $\|E_N\|_\infty$ uniformly bounded by $||E||_\infty$,
any  Hermitian matrices $(D_N,\tilde D_N)_{N\in\N}$,
$$\left|{1\over N}\log I_N^{(\beta)}(D_N,E_N)
-{1\over N}\log I_N^{(\beta)}(\tilde D_N,E_N)
\right|\ppq ||E||_\infty \tr |D_N-\tilde D_N|$$
\item 
Let $D_N=\mbox{diag }(\theta,0,\cdots, 0)$.
Assume that 
there is a positive $\eta$
and a finite integer $N_0$ such that for $N\ge N_0$, 
$\frac{2 \theta}{\b}\in H_{E_N}( [\l_{\rm min}(E_N)-
\eta,\l_{\rm max}(E_N)+\eta]^c)$.
We
let $v_N$ be the unique solution in $-\b(2\theta)^{-1}+
 [\l_{\rm min}(E_N)-
\eta , 
\l_{\rm max}(E_N)+\eta]^c$ of the equation 
\begin{equation}
\label{defvN}
{\b \over 2\theta} H_{E_N}\left(\frac{\b}{2\t}+v_N \right)=1.
\end{equation}

Then,  $v_N \in [\l_{\rm min}(E_N), \l_{\rm max}(E_N)]$ and 
for any $\zeta\in (0,{1\over 2})$,
there exists a finite constant $C(\eta,\zeta)$
depending only on $\eta$ and $\zeta$ such that 
for all $N\ge N_0$
$$
\left|{1\over N} \log I_N^{(\beta)}(\theta,E_N)-\theta  v_N +{\beta\over 2N}
\sum_{i=1}^N \log\left(1+{2\t\over \beta}
 v_N -{2\t\over \beta}
 \l_i\right) \right|
\ppq C(\eta,\zeta)N^{-{1\over 2} +\zeta} \|E_N\|_\infty. $$

\item Let $D_N=\mbox{diag}(\theta,0,\cdots, 0)$.
Let $E_N,\tilde E_N$ be two matrices such
that 
$$d(\mun_{E_N},\mun_{\tilde E_N})\ppq \d,$$
where $d$ is the Dudley distance on $\PR$
and so that both $E_N$ and $ \tilde E_N$ satisfy \eqref{hypnorm}.\\
Let $\eta>0$.
Assume that there exists $N_0<\infty$ so that for $N\ge N_0$,
$\frac{2 \theta}{\b}\in H_{E_N}( [\l_{\rm min}(E_N)-
\eta,\l_{\rm max}(E_N)+\eta]^c)
\cap  H_{\tilde E_N}( [\l_{\rm min}(\tilde E_N)
-\eta,\l_{\rm max}(\tilde
E_N)+\eta]^c)$.
Then, there exists a function $g(\d,\eta)$ (independent of $N$)
going to zero with $\d$ for any $\eta$ and 
such that for all $N\ge N_0$
$$\left|{1\over N}\log I_N^{(\beta)}(D_N,E_N)
-{1\over N}\log I_N^{(\beta)}(D_N,\tilde E_N)
\right|\ppq g(\d,\eta)$$
\end{enumerate}

\end{lem}
Note that
the third point is analogous to the continuity
statement
obtained in the case where $D_N$ 
has also rank $N$ in \cite{GZ}, Lemma 5.1.
However, let us mention again that
there is an important difference here 
which lies in the fact that the smallest
and largest eigenvalues play quite an important role.
In fact, it can be seen (see Theorem \ref{expreel})
that if we let one eigenvalue
be much larger than the support of
the limiting spectral distribution,
then the
limit of the spherical integral 
will change dramatically. 
However, Lemma \ref{cont}.3 shows that this
limit will not depend on these escaping eigenvalues 
provided $|\theta|$ is smaller than some critical value
$ \theta_0(\l_{\rm min},\l_{\rm max})$ ($= \min(|H_{min}\b/2|, |H_{max}\b/2|)$). \\
\indent Before going into the proof of Lemma \ref{cont}, 
let us show that Theorem \ref{sh} is a direct consequence 
of its second point.

\nn
{\bf Proof of Theorem \ref{sh} :}
Since we assumed that, for $N$ large enough,
 $2\theta \b^{-1}\in H_{E_N}([\l_{\rm min}(E_N)-\eta,\l_{\rm max}(E_N)+\eta]^c) ,$
we can find a $v_N$ satisfying (\ref{defvN}). Note that $v_N$ is unique by strict 
monotonicity of $ H_{E_N}$ on $]-\infty \, , \, \l_{\rm min}(E_N)
 -\eta[$,
where it is negative, and on $]\l_{\rm max}(E_N)+\eta \, , \,
\infty[$, 
where it is positive.  Therefore, 
$$ (2\theta)^{-1}+v_N\in [\l_{\rm min}(E_N)-\eta \, , \, \l_{\rm max}(E_N)+\eta]^c
$$
 ensures that 
\begin{equation}\label{boundv}
1-\frac{2\theta}{\b}\l_i+ \frac{2\theta}{\b} v_N> \frac{2|\theta|}{\b}\eta
\end{equation}
so that, 
because
of the uniform continuity of  $H_{E_N}$ 
on $[\l_{\rm min}( E_N)
-\eta,\l_{\rm max}(
E_N)+\eta]^c$, as $\mun_{E_N}$ converges to $\mu_E$, 
$v_N$ converges to $v$ the solution of 
$\HE\left( \frac{\b}{2\t} +v\right) = \frac{2\t}{\b}$
and
$$\lim_{N\ra\infty}{1\over N}\sum_{i=1}^N\log\left(  
1+{2\t \over \beta}
 v_N -{2\t\over \beta}
 \l_i\right)=\int \log\left( 1+{2\t\over \beta}
 v -{2\t\over \beta}
 \l\right) d\mu_E(\l).$$
Furthermore, the computation
of the derivative of
$$\t \mapsto \t v - \frac \b 2 \int \log\left( 1+ \frac{2\t}{\b}v - \frac{2\t}{\b}\l\right)
d\mu_E(\l),$$
with this particular $v=\RE(2\theta\beta^{-1}
)$ allows us to get the explicit expression 
$$\t v - \frac \b 2 \int \log\left( 1+ \frac{2\t}{\b}v - \frac{2\t}{\b}\l\right)
d\mu_E(\l)=\frac \b 2\int_0^{2\t\over\b} \RE(u)du.$$
Therefore, Hypothesis \eqref{hypnorm}
together with Lemma \ref{cont}.2 finishes the proof of 
(\ref{int1}).\\

Now the last point is to check that 
under Hypothesis \ref{hypsurEN}, 
the assumption of  Lemma \ref{cont}.2 is equivalent to $2\t/\beta\in
\HE([\l_{\rm min},\l_{\rm max}]^c)$.\\
Let us first observe that 
$\HE([\l_{\rm min},\l_{\rm max}]^c)
=\bigcup_{\eta>0} \HE([\l_{\rm min}-\eta,\l_{\rm max}+\eta]^c)$
and that, under Hypothesis \ref{hypsurEN},
$$\HE([\l_{\rm min}-2\eta,\l_{\rm max}+2\eta]^c)\subset
\bigcup_{N_0\ge 0} \bigcap_{N\ge N_0}
H_{E_N}([\l_{\rm min}(E_N)-\eta,\l_{\rm max}(E_N)+\eta]^c),$$
since, for any $\l \in \bigcup_{N_0\ge 0} \bigcap_{N\ge N_0}
\mbox{supp}(\mun_{E_N})$, the application $z \mapsto (z-\l)^{-1}$
is continuous bounded on $[\l_{\rm min}-2\eta,\l_{\rm max}+2\eta]^c$.
Therefore, $\frac{2\t}{\beta}\in\!
\HE([\l_{\rm min},\l_{\rm max}]^c)$
implies the assumption of  Lemma \ref{cont}.2.\\
Conversely, we get by the same arguments that
$$
\bigcup_{N_0\ge 0} \bigcap_{N\ge N_0}
H_{E_N}([\l_{\rm min}(E_N)-2\eta,\l_{\rm max}(E_N)+2\eta]^c)\subset
\HE([\l_{\rm min}-\eta,\l_{\rm max}+\eta]^c),$$
what completes the proof.
\hfill\xx

\subsection{ Proof of Lemma \ref{cont}}

\noindent $\bullet$ The first point is trivial since the matrix $U$ is 
unitary or orthogonal and hence bounded.\\
$\bullet$ Let us consider the second
point. %We can assume without
%loss of generality that $E_N$
%is  non-negative (we translate all the 
%$\l_i$'s by the same constant if necessary). 
We now
stick to the case $\b=1$
and will summarize
at the end of the proof
the changes to perform
for the case $\b=2$.
We can assume that the $\{\l_1(E_N),\cdots,\l_N(E_N)\}$
is not reduced to
a single point $\{e\}$ since otherwise
the result is straightforward. We write
in short $I_N(\theta,E_N)=I_N^{(1)}(D_N, E_N)$.
The ideas of the proof are very close 
to usual large deviations
techniques, and in fact in some sense simpler
because strong concentration  arguments 
are available for free (cf. (\ref{rho})).
 Following Fact \ref{gaussien},
we can write, with $(\l_1,\cdots,\l_N)$ the eigenvalues
of $E_N$, 
$$I_N(\theta,E_N)=\E\left[ \exp\left\{N\theta{\sum_{i=1}^N \l_i
g_i^2\over \sum_{i=1}^N 
g_i^2}\right\} \right]$$
where the $g_i$'s are i.i.d standard Gaussian variables.
Now, writing the Gaussian  vector $(g_1, \ldots, g_N)$
in its polar decomposition, we realize 
of course that the spherical integral does not depend on 
its radius $r=\|g\|$ which follows the law 
$$\rho_N(dr):=Z_N^{-1}r^{N-1}e^{-{1\over 2} r^2}dr,$$
with $Z_N$ the appropriate normalizing constant.\\
The idea of the proof is now that $r$ will of course concentrate around $\sqrt N$ 
so that we are reduced to study the numerator and to make the adequate 
change of variable so that it concentrates around $v_N$.
For $\kappa<1/2$, there exists a finite constant
$C(\kappa)$ such
that
\begin{equation}
\label{rho}
\rho_N\left(\left|{r^2\over N}-1\right|\pgq N^{-\kappa}\right)\ppq
C(\kappa) e^{ -{1\over 4}N^{1-2\kappa}}.
\end{equation}
Such an estimate can be readily obtained 
by applying standard precise Laplace method to the
law $\tilde\rho_N$ of $(N-2)^{-1} r^2$ which is given by
$$\tilde\rho_N(dx)=\tilde Z_N^{-1} 1_{x\ge 0} e^{-\frac{N-2}{2} f(x)} dx$$
with $f(x)= x -\log x$. Indeed, $f$ achieves its minimal value at $x=1$
so that for any $\e>0$, there exists $c(\e)>0$
such that  $\tilde Z_N\tilde\rho_N(|x-1|>\e)\le e^{-c(\e)N}$.
Now, $\s_\e=\inf\{f"(x), |x-1|\le \e\}>0$
so that Taylor expansion results with
$$\tilde Z_N\tilde\rho_N(|x-1|\ge N^{-\kappa})\le e^{-c(\e)N} +
\int_{y\ge  N^{-\kappa}}
e^{-\frac{N-2}{2}\s_\e y^2}dy\le e^{-\frac{\s_\e}{3} N^{1-2\kappa}}$$
where the last inequality holds for $N$ large enough.
A lower bound on $\tilde Z_N$ is obtained similarly by considering 
$\tilde\s_\e= \sup\{f"(x), |x-1|\le \e\}>0$
showing that $\tilde Z_N\ge \tilde c(\e)\sqrt{N}^{-1}$. We conclude by noticing that $\s_\e$
goes to one as $\e$ goes to zero. Note that such a result can also 
be seen as a direct consequence of moderate
deviations (cf. section 3.7 in \cite{DZ}). %Ader

From this, if we introduce the event 
$A_N(\k) := \left\{ \left|{\|g\|^2\over N}-1\right|\ppq N^{-\kappa}\right\},$
it is not hard to see that for
any $\k < \frac 1 2$ and for $N$ large enough (such that
$1-C(\kappa) e^{ -{1\over 4}N^{1-2\kappa}} >0$), we have
\begin{equation}\label{conco}\nonumber
 1\ppq{ I_N(\theta,E_N)\over
\E\left[
1_{A_N(\k)}
\exp\left\{N\theta
{\sum_{i=1}^N \l_i g_i^2\over \sum_{i=1}^N  g_i^2}
\right\}\right]}
\ppq\d(\kappa,N)
\end{equation}
where 
$\d(\kappa,N)={1\over 1- C(\kappa) e^{ -{1\over 4}N^{1-2\kappa}}}$.  
Therefore,
\begin{eqnarray}
\!\!\!\!\!\!\!\!\!\!\!\!\!I_N(\theta,E_N) & \ppq & \d(\kappa,N)\E\left[
1_{A_N(\k)}
\exp\left\{N\theta
{\sum_{i=1}^N \l_i g_i^2\over \sum_{i=1}^N  g_i^2}
\right\}\right]  \label{po1} \\
& \ppq & \d(\kappa,N)e^{N \theta  v +N^{1-\kappa}|\theta |(\|E_N\|_\infty+|v|) } 
\E\left[1_{A_N(\k)}
 \exp\left\{ \theta {\sum_{i=1}^N \l_i
g_i^2}-v \theta \sum_{i=1}^N 
g_i^2\right\}\right]
\nonumber
\end{eqnarray}
for any $v\in\RR$.
Now,
\begin{equation}
\E\left[1_{A_N(\k)}
 \exp\left\{\theta {\sum_{i=1}^N \l_i
g_i^2}-v \theta \sum_{i=1}^N 
g_i^2\right\}\right] = \prod_{i=1}^N\left[\sqrt{ 1+2 \theta v -2 \theta\l_i}\right]^{-1} \,\,
 P_N(A_N(\k))\label{po2}
\end{equation}
with $P_N$ the
probability measure on $\RR^N$ given by 
$$P_N(dg_1, \ldots,dg_N)=\frac{1}{\sqrt{2\pi}^N}
\prod_{i=1}^N \left[\sqrt{ 1+2 \theta v -2 \theta\l_i}\,\,
e^{-{1\over 2}(  1+2 \theta v -2 \theta\l_i)g_i^2}dg_i \right]$$
which is well defined provided 
we choose $v$ so that 
\begin{equation}\label{cond1}
 1+2 \theta v -2 \theta\l_i>0\quad\forall \,i \mbox{ from }  1 \mbox{ to } N.
\end{equation}
Thus, for any such $v$'s, we get from (\ref{po1}) and (\ref{po2}),
that for any $\k = \frac 1 2 - \zeta$ with $\zeta >0$
and $N$ large enough, since $P_N(A_N(\k))\le 1$,
\begin{equation}
I_N(\theta,E_N) \ppq \d(\kappa,N) \prod_{i=1}^N \left[\sqrt{ 1+2 \theta v -2 \theta\l_i}\right]^{-1}
 e^{N \theta  v+N^{1-\k}|\t v|}
e^{N^{1-\k} |\t| \|E_N\|_\infty}.
\label{re1}
\end{equation}
We similarly 
obtain
the lower bound
 \begin{equation}
\nonumber
I_N(\theta,E_N)\pgq
e^{N \theta  v - N^{1-\kappa}|\theta | (\|E_N\|_\infty +|v|)}
\prod_{i=1}^N\left[\sqrt{ 1+2 \theta v -2 \theta\l_i}\right]^{-1}P_N(A_N(\k))
\label{re2}
\end{equation}
Now, we show that we
can choose $v$ wisely so that for $N\pgq N(\kappa)$,
\begin{equation}\label{desired}
P_N(A_N(\k)) = P_N(|N^{-1}||g||^2-1|\ppq N^{-\kappa}) \pgq {1\over 2}.
\end{equation}
This will  finish to prove, with this choice of $v$,
that
\begin{equation}
I_N(\theta,E_N)
\pgq {1\over 2}
e^{N \theta  v - N^{1-\kappa}|\theta | (\|E_N\|_\infty +|v|)}
\prod_{i=1}^N\left[\sqrt{ 1+2 \theta v -2 \theta\l_i}\right]^{-1}
\label{re3}
\end{equation}
yielding the desired  lower bound.\\
We know that
$P_N$ is a product measure
under which
$$\tilde g_i=\sqrt{1+2 \theta v -2 \theta\l_i}\,\,g_i$$
are i.i.d standard Gaussian variables.
Let us now choose $v=v_N$  in 
$-(2\theta)^{-1}+ [\l_{\rm min}(E_N)-\eta \, , \, \l_{\rm max}(E_N)+\eta]^c$
satisfying
\begin{equation}
\label{veq}
\E_{P_N}\left[{1\over N}\|g\|^2\right]=
\E\left[{1\over N}\sum_{i=1}^N  {\tilde g_i^2\over 1+2 \theta v_N -2 \theta\l_i}\right]
={1\over 2\theta} H_{E_N}\left( (2\theta)^{-1}+v_N\right)=1.
\end{equation}
We recall from \eqref{boundv} that $
1-2\theta\l_i+2\theta v_N>2|\t|\eta>0$
so that all our computations are validated by this final choice.

With this choice of $v_N$, we have
$$\E_{P_N}\left[\left({1\over N}\|g\|^2-1\right)^2\right]
={2\over N^2}\sum_{i=1}^N {1\over (1+2 \theta v_N -2 \theta\l_i)^2}\ppq {2\over N \t^2\eta^2}$$
so that by Chebychev's inequality
$$P_N(|N^{-1}||g||^2-1|\pgq N^{-\kappa})\ppq {2\over \eta^2 \t^2}N^{2\kappa-1},$$
which is smaller than $2^{-1}$ for sufficiently large $N$ since $2\kappa<1$,
resulting with (\ref{desired}).\\
Finally, since by definition
$${1\over N}\sum_{i=1}^N
\frac{1}{ 1-2\theta \l_i+2\theta v_N}
=1$$
with $(\l_i)_{1\le i\le N}$ which do not all take the
same value, there exists
$i$ and $j$ so that
$$ -2\theta \l_i+2\theta v_N>0,\quad -2\theta \l_j+2\theta v_N<0$$
so that $v_N\in [ \l_{\rm min}(E_N) ,\l_{\rm max}(E_N)]$.
Thus, (\ref{re3}) together with (\ref{re1}) give the second point of the lemma for $\beta=1$.\\

In the case where $\beta=2$, the $g_i^2$
have to be replaced everywhere
by $g_i^2+\hat g_i^2$
with independent Gaussian variables $(g_i,\hat g_i)_{1\ppq i\ppq N}$.
This time, we can concentrate
$$\frac 1 N \|g\|^2={1\over N}\sum_{i=1}^N g_i^2
+{1\over N}\sum_{i=1}^N \hat g_i^2$$
around $2$. Everything then follows by dividing
$\t$ by two
and noticing that we will get the
same Gaussian integrals squared.\\

$\bullet$ The
 last point is an  easy consequence
of the second  since, 
for any $\l \in \bigcup_{N_0\ge 0} \bigcap_{N\ge N_0}
(\mbox{supp}(\mun_{E_N})$ $\cap \,\, \mbox{supp}(\mun_{\tilde E_N}))$, the application $z \mapsto (z-\l)^{-1}$
is continuous bounded (with norm depending on $\eta$)
 on $\bigcup_{N_0\ge 0} \bigcap_{N\ge N_0} [\l_{\rm min}(E_N)-\eta,\l_{\rm max}(E_N)+\eta]^c$. \hfill\xx

\subsection{Generalisation of the method to the multi-dimensional case}
In the sequel, we want to apply the strategy we used above to show
Theorem \ref{rgfini}, that is to say study the behaviour of the spherical integrals
as the rank of $D_N$ remains negligible compared to $\sqrt N$. In this case
and if all the eigenvalues of $D_N$ are small enough, we show that
it behaves like a product, namely that we have the equality \eqref{sip}.  
To lighten the notations, we let $\t_i := \t_i^N$, for all $i \ppq M(N)$.\\
We will rely again on Fact \ref{gaussien} and write
in the case $\b =1$,
\begin{equation}
\label{gmd}
I_N(D_N, E_N) = \EE\left[ \exp \left\{ N \sum_{m=1}^M \t_m \frac{\sum_{i=1}^{N} \lambda_i (\tilde g_i^{(m)})^2}{\sum_{i=1}^{N} (\tilde g_i^{(m)})^2} \right\}\right],
\end{equation}
where the expectation is taken under the standard Gaussian measure and the vectors $(\tilde g^{(1}), \ldots, \tilde g^{(M)})$
are obtained from the Gaussian vectors $(g^{(1)}, \ldots, g^{(M)})$ by a standard Schmidt orthogonalisation procedure.\\
This means that there exists a lower triangular matrix $A = (A_{ij})_{1 \ppq i,j \ppq M}$ such that
for any integer $m$ between $1$ and $M$, 
$$ \tilde g^{(m)} = g^{(m)} + \sum_{j=1}^{m-1} A_{mj} g^{(j)} $$
and the $A_{ij}$'s are solutions of the following system :
for all $p$ from $1$ to $m-1$,
\begin{equation}
\label{syst}
\langle g^{(m)}, g^{(p)} \rangle + \sum_{j=1}^{m-1} A_{mj} \langle g^{(j)}, g^{(p)} \rangle = 0,
\end{equation}
with $\langle.,.\rangle$ the usual scalar product in $\RR^N$.\\

Therefore, if we denote, for $i$ and $j$ between $1$ and $M$, with $i \ppq j$,
$$ X_N^{ij} := \frac 1 N \langle g^{(i)}, g^{(j)}\rangle$$
and
$$ Y_N^{ij} := \frac 1 N  \sum_{l=1}^N \l_l g^{(i)}_l g^{(j)}_l,$$
then, for each $m$ from $1$ to $M$, there exists a rational function $F_m : \RR^{m(m+1)} \ra \RR$
such that 
\begin{equation}
\label{eqFk}
\frac{\sum_{i=1}^N \lambda_i (\tilde g_i^{(m)})^2}{\sum_{i=1}^{N} (\tilde g_i^{(m)})^2} 
= F_m(( X_N^{ij},  Y_N^{ij})_{1 \ppq i \ppq j \ppq m})
\end{equation}
  and a rational function $G_m : \RR^{\frac{m(m+1)}{2}} \ra \RR$ such that
\begin{equation}
\label{eqGk}
\frac 1  N \sum_{i=1}^{N} (\tilde g_i^{(m)})^2  
= G_m(( X_N^{ij})_{1 \ppq i \ppq j \ppq m}).
\end{equation}

We now adopt the following system of coordinates
in $\RR^{MN}$ : $r_1$, $\a^{(1)}_1, \ldots, \a^{(1)}_{N-1}$
are the polar coordinates of $g^{(1)}$,
$r_2 := \|g^{(2)}\|$, $\b_2$ is the angle between
$g^{(1)}$ and $g^{(2)}$,  
$\a^{(2)}_1, \ldots, \a^{(2)}_{N-2}$ are the angles needed to 
spot $g^{(2)}$ on the cone of angle $\b_2$ around
$g^{(1)}$,
then  $r_3 := \|g^{(3)}\|$, $\b_3^i$ the angle between
$g^{(3)}$ and $g^{(i)}$ ($i= 1, \, 2$) and
$\a^{(3)}_1, \ldots, \a^{(3)}_{N-3}$ the angles needed to
spot $g^{(3)}$ on the intersection of the two cones...etc...\\
Then observe that  $F_m(( X_N^{ij},  Y_N^{ij})_{1 \ppq i \ppq j \ppq m})$
depends only on the $\a$'s (because the 
$\frac{\tilde g^{(i)}}{\|\tilde g^{(i)}\|}$ do) whereas
$G_m(( X_N^{ij})_{1 \ppq i \ppq j \ppq m})$ depends on the
$r$'s and the $\b$'s. Therefore, if we consider the event
$$ B_{N}(\k) := \left\{ \forall i, \quad \left|X_N^{ii} -1\right| \ppq N^{-\k}, \quad \forall i \neq j
\quad \left| X_N^{ij}\right|\ppq N^{-\k}
\right\},$$
then, as in the case of rank one, we can write that
\begin{equation}
\label{dec}
I_N(D_N, E_N) \ppq \EE\left[1_{B_N(\k)}
 e^{N \t_m F_m(X_N^{ij},  Y_N^{ij})}\right] + P(B_N(\k)^c) I_N(D_N, E_N).
\end{equation}
Now we claim that, for $N$ large enough, for any $\k > 0$,
 there exists an $\a > 0$ such that
\begin{equation}
\label{devk}
P(B_N(\k)^c) \ppq C^\prime(\k) e^{-\a N^{1-2\k}}.
\end{equation}
Indeed, as in (\ref{rho}),
\begin{eqnarray*}
P(B_N(\k)^c) & \ppq & \sum_{i=1}^M P\left(\left| X_N^{ii} - 1\right| > N^{-k} \right) + \sum_{i, j=1}^M P\left(\left| X_N^{ij}\right| > N^{-k} \right) \\ 
& \ppq & c_1(\k) M e^{-\frac 1 4 N^{1-2\k}} + c_2(\k) M^2 e^{-\frac 1 2 N^{1-2\k}},
\end{eqnarray*}
what gives immediately (\ref{devk}).\\
Now, as far as $\k < \frac 1 2$, (\ref{dec}) together
with (\ref{devk}) give 
$$ 1 \ppq \frac{I_N(D_N, E_N)}{\EE\left[1_{B_N(\k)}
 e^{N F_m(X_N^{ij},  Y_N^{ij})}\right]}\ppq 1+\e(N,k),$$
with $\e(N,k)$ going to zero.

We now want to expand $F_M$ on $B_N(\k)$ as we did 
in the previous subsection. \\
As the $A_{ij}$'s satisfy the linear system
 (\ref{syst}), we can write the Cramer's formulas
corresponding to it and get
$$ A_{ij} = \frac{\mbox{det}(R_N^{kl})_{1 \ppq k,l \ppq i-1} }
{\mbox{det}(X_N^{kl})_{1 \ppq k,l \ppq i-1} },$$
where 
\[
R_N^{kl} = \left\{ \begin{array}{ll}
X_N^{kl}, & \textrm{ if } l \neq j \\
- X_N^{ki} & \textrm{ if } l = j.
\end{array}
\right.
\]
Now, we look at the denominator and can show that
$$ \mbox{det}(X_N^{kl})_{1 \ppq k,l \ppq i-1}
\pgq 1 - \sum_{s = 1 }^{i-1} (MN^{-\k})^s \pgq \frac 1 2, $$
where the last inequality holds for
$N$ large enough as far as $M = o(N^{\k})$.\\
We now go to the numerator : expanding over the
$j$th column, we get this time that
$$ \mbox{det}(R_N^{kl})_{1 \ppq k,l \ppq i-1}
\ppq N^{-\k} +  (M-1) N^{-2\k} \sum_{s = 1 }^{i-1} (MN^{-\k})^s
\ppq c N^{-\k},$$
where again the last equality holds as far as $M = o(N^{\k})$
and $c$ is a fixed constant.\\
From the two last inequalities, we have that, on $B_N(\k)$,
$\sup_{i<j} |A_{ij}| 
\ppq c^{\prime} N^{-\k}.$ \\
From that we can easily deduce that, for any
$m$ less than $M$, we have \vspace{-0.3cm}
\begin{eqnarray*}
  \frac 1 N \left\| \tilde g^{(m)} -  g^{(m)}\right\|^2
& \ppq & \frac 1 N \sum_{i,j=1}^{m-1} |A_{mj} A_{mi}| |\langle g^{(i)}, g^{(j)}\rangle|^2 \\
& \ppq & c^{\prime\prime} N^{-2\k} (M^2 N^{-2\k} + M) 
\ppq c_3 N^{-\k} . 
\end{eqnarray*}

From these estimations and (\ref{gmd}), for any $v_j^N$, we get the following upper bound~:
\begin{multline*}
I_N(D_N, E_N)  \ppq  (1+\e(\k, N))  \exp\Bigl\{N \sum_{j=1}^M \t_j v_j^N\Bigr\}\\
\EE\left[1_{B_N(\k)} \prod_{j=1}^M 
\exp\left\{ N\t_j \frac{ \frac 1 N\sum_{i=1}^N \l_i \bigl(\tilde g_i^{(j)}\bigr)^2 
- v_j^N \frac 1 N \sum_{i=1}^N  \bigl(\tilde g_i^{(j)}\bigr)^2}{1 
+\frac 1 N \left( \|\tilde g^{(j)}\|^2 - \|g^{(j)}\|^2\right) +
\left( \frac 1 N  \|g^{(j)}\|^2 - 1 \right) } \right\}
\right] 
\end{multline*}
\vspace{-0.3cm}
\begin{multline*}
\ppq  (1+\e(\k, N)) \exp\Bigl\{N \sum_{j=1}^M \t_j v_j^N\Bigr\}\\
\EE\left[1_{B_N(\k)} \prod_{j=1}^M 
\exp\left\{ \left(\t_j \sum_{i=1}^N \l_i \bigl(\tilde g_i^{(j)}\bigr)^2 
- v_j^N \t_j \sum_{i=1}^N  \bigl(\tilde g_i^{(j)}\bigr)^2\right) \left[ 1 + c_4  N^{-\k} \right] 
\right\}
\right]
\end{multline*}
\vspace{-0.3cm}
\begin{multline*}
\ppq (1+\e(\k, N))\, e^{\bigl\{N \sum_{j=1}^M \t_j v_j^N\bigr\}}\, 
e^{\bigl\{ C \sup|\t_j| (\|E_N\|_\infty + \sup|v_j^N|) M N^{1-\k}\bigr\}}\\
\EE\left[ \prod_{j=1}^M 
\exp\left\{ \t_j   \sum_{i=1}^N \l_i \bigl( g_i^{(j)}\bigr)^2 
- v_j^N \sum_{i=1}^N  \bigl( g_i^{(j)}\bigr)^2
\right\},
\right].
\end{multline*}
where $C$ is again a fixed constant.\\
From the hypotheses of Theorem \ref{rgfini},
we know that there exists an $N$ such that
$2\t_j \in H_{E_N}([\l_{\rm min}(E_N)-\eta, \l_{\rm max}(E_N)+\eta]^c), $
from which we can easily deduce that $|2\t_j | \ppq \eta^{-1}$.
Moreover, as in the proof of Lemma \ref{cont}.2,
$|v_j^N|\le ||E_N||_\infty$ 
is  uniformly bounded. Therefore, we get
$$ \limsup_{N \ra \infty} {1\over NM(N)} \log  I_N(D_N, E_N)  \ppq \int \IE(\t) d\mu_D(\t). $$
We also get a similar lower bound and conclude similarly
to the preceding subsection
by considering
the shifted probability measure
$P_N^{\t_1, \ldots,\t_M} = \otimes_{j=1}^M P_N^{\t_j}$ where
$$ P_N^{\t_j}(dg_1, \ldots, dg_N) = \frac{1}{\sqrt{2\pi}^N}
\prod_{i=1}^N \sqrt{1+2\t_j v_j^N -2\t_j\l_i}\,\,
e^{-\frac 1 2 (1+2\t_j v_j^N -2\t_j\l_i) g_i^2} dg_i.$$
This concludes the proof of Theorem \ref{rgfini}. \hfill \xx

%%%%%%%%%%%%%%%%%%%%%%%%%%%%%%%%%%%%%%%%%%%%%%%%%%%%%%%%%%%%%%%%%%%%%%%%%%%%%%%%%%%%%%%%%%

\section{Central limit theorem in the case of rank one}
\label{tcl}

Under the hypotheses of Theorem \ref{sh}, $v_N$ (defined by (\ref{defvN})) 
is converging to $v = \RE\Bigl(\frac{2\t}{\b}\Bigr)$ and we established 
that the spherical integral is converging to
$\theta  v -{\beta\over 2}
\int \log\bigl(1+{2\t\over \beta}
 v -{2\t\over \beta}
 \l\bigr) d\mu_E(\l)$.
 In the case where the fluctuations of the eigenvalues
do not interfere, we can get sharper estimates,
given, in the case $\b=1$, by Theorem \ref{precise1}.
This section is devoted to its proof, namely
the study of the behaviour of
$e^{-N\bigl(\t  \RE(2\t) - \frac{1}{2N} \sum \log(1+2\t \RE(2\t) - 2\t \l_i) \bigr)} I_N(\t, E_N)$.\\

\nn
\textbf{Proof of Theorem \ref{precise1}}

$\bullet$
We first treat the
non degenerate case $\mu_E\neq \d_e$.

Let us first make an important remark : the hypothesis
that $d(\mun_{E_N}, \mu_E) = o(\sqrt N^{-1})$ has the two 
following consequences :
\begin{equation}
\label{vproche}
|v-v_N| = o(\sqrt N^{-1})
\end{equation}
\begin{equation}
\label{pesah}
\textrm{and } \quad \lim_{N \ra \infty} \sqrt N (H_{E_N} - \HE)(\KE(2\t)) =0.
\end{equation}
Indeed, since $2\t \in \HE([\lm, \lM]^c)$,
there is an $\eta >0$, such that, for $N$ large enough,
 $2\t \in H_{E_N}([\l_{\rm min}(E_N) - \eta, \l_{\rm max}(E_N)+\eta]^c). $
Therefore, as for any $\l$ which is in $\mbox{supp}(\mun_{E_N})$
for $N$ large enough,
$z \mapsto (z-\l)^{-1}$ is uniformly bounded and Lipschitz
on $\bigcap_{N \pgq N_0} [\l_{\rm min}(E_N)-\eta, \l_{\rm max}(E_N)+\eta]^c,$
we get directly \eqref{vproche}, and also
\eqref{pesah} as we know that $\KE(2\t) \in [\lm, \lM]^c.$\\
 
\nn
For $v=\RE(2\theta)$, we set
$$ \g_N=\left( {1\over N}\sum_{i=1}^N g_i^2-1\right) \textrm{ and }
\hat \g_N=\left( {1\over N}\sum_{i=1}^N \l_i
g_i^2-v\right).$$
Let us also define
for $\e>0$
$$I_N^{\e}(\t,E_N):=\int_{|\g_N|\ppq\e, |\hat \g_N|\ppq\e}
\exp\left\{ \theta N {\hat \g_N + v\over \g_N+1}\right\}
\prod_{i=1}^N dP(g_i),$$
with $P$ the standard Gaussian probability measure on $\RR$.
We claim that,
for any $\zeta>0$, for $N$ large enough,
\begin{equation}\label{contr}
\left|I_N(\t,E_N)
- I_N^{\e}(\t,E_N) \right|\ppq e^{- N^{1-\zeta}} I_N(\t,E_N).
\end{equation}
Indeed,  consider
$$\mu_N^\t (dg)={1\over I_N(\theta,E_N)}
\exp\left\{ \theta N {\sum_{i=1}^N \l_i g_i^2
\over \sum_{i=1}^N  g_i^2}\right\}
\prod_{i=1}^N dP(g_i).$$
(\ref{contr}) is equivalent
to
\begin{equation}\label{contr2}
\mu_N^\t\left(|\g_N|\ge\e\right)\le {1\over 2} e^{- N^{1-\zeta}}
\mbox{ and } 
\mu_N^\t\left(|\hat\g_N|\ge\e\right)\le {1\over 2} e^{-
N^{1-\zeta}}
\end{equation}
The first inequality is trivial since by (\ref{rho}), for
$\kappa<\frac 12$,
$$\mu_N^\t\left(|\g_N|\ge N^{-\kappa}\right)=
\rho_N\left(\left|\frac {r^2}{ N} -1\right|\ge N^{-\kappa}\right)
\le  e^{-\frac 14 N^{1-2\kappa}}.$$
To show the second point, following the proof of Lemma
\ref{cont}, we find a finite
constant $C(\kappa)$   so that
$$\mu_N^\t\left(|\hat\g_N|\ge\e\right)
\le C(\kappa) e^{C(\kappa) N^{1-\kappa}|\t|||E_N||_\infty}
P_N\left( |\hat\g_N|\ge\e\right)$$
where under $P_N$ the $g_i$ are independent centered 
Gaussian variable with covariance $(1-2\theta \l_i+2\theta v_N)^{-1}$.
Hence
$$P_N\left( |\hat\g_N|\ge\e\right)
=P^{\ot N}\left( \left| {1\over N}\sum_{i=1}^N \frac{\l_i}{
 1-2\theta \l_i+2\theta v_N} \tilde g_i^2 -v\right|\ge \e\right).$$
Let us denote $\tilde E_N=\phi_{v_N}( E_N)
$ with $\phi_v(x)= x(1-2\t x+2\t v)^{-1}$.
Then, 
the spectral measure of $\tilde E_N$
converges towards $\mu_{\tilde E} :=\phi_v\sharp \mu_E$ since
$v_N$ converges
towards $v$ (see (\ref{vproche})).
Moreover $\l_{\rm min}(\tilde E_N)$
and $\l_{\rm max}(\tilde E_N)$
converge.
Hence, we can apply Lemma
\ref{lem2} 
to obtain a large deviation principle
for the law of ${1\over N}\sum_{i=1}^N \l_i(\tilde E_N)\tilde g_i^2$ 
under $P^{\ot N}$ with good rate function $J$.
One checks that
$J$ has a unique minimizer which is 

$$z_0=R_{\mu_{\tilde E}}(0)=
 \int {\l\over 1-2\theta \l +2\theta v}d\mu_E(\l)= v.$$
As a consequence, for $\e>0$, there exists $\d(\e)>0$ so that
for $N$ large enough 
$$P^{\ot N}\left( \left| {1\over N}\sum_{i=1}^N \frac{\l_i}{
 1-2\theta \l_i+2\theta v_N} \tilde g_i^2 -v\right|>\e\right)
\le e^{-\d(\e) N}.$$
This completes the proof of (\ref{contr2}).

\bigskip

We now deal with $I_N^{\e}(\t,E_N)$.
We use the expansion
 $ \frac{1}{1+\g_N} = 1-\g_N + \frac{\g_N^2}{1+\g_N}$ 
to get that
$$
I_N^{\e}(\t,E_N)=e^{N\theta v}
\int_{|\g_N|\ppq\e, |\hat \g_N|\ppq\e}
e^{\left\{ -\theta N \g_N {\hat \g_N - v\g_N\over \g_N+1}\right\}}
\,e^{\{\theta N(\hat \g_N- v\g_N)\}} \prod_{i=1}^N dP(g_i).
$$
We note that
$$\exp\{\theta N(\hat \g_N- v\g_N)\} \prod_{i=1}^N dP (g_i)
=\prod_{i=1}^N \left[\sqrt{1+2\theta v-2\t\l_i}\right]^{-1} \prod_{i=1}^N dP_i(g_i)$$
with $P_i$ the centered Gaussian probability measure
$$dP_i(x)=\sqrt{(2\pi)^{-1}(1+2\theta v-2\t \l_i)}\,\,
\exp\left\{ -{1\over 2} (1+2\theta v-2\t \l_i) x^2\right\}dx.$$
We have that
\begin{equation}\label{eqp}
1+2\theta v-2 \t \l_i=2\theta(\KE(2\theta)-\l_i)
\end{equation}
and we know that $\KE(2\theta)\in [\l_{\min},\l_{\rm max}]^c
$ .
Further, arguing as in (\ref{boundv}),
 we  find, for any given $\theta>0$,
a constant $\eta_\theta>0$ such that
$$\inf_{1\ppq i\ppq N}(1+2\theta v-2 \t \l_i)>\eta_\theta$$
insuring that the $P_i$ are well defined.
Therefore, 
\begin{multline}
\label{presque}
I_N^{\e}(\t,E_N) = e^{N\theta v-{N\over 2}
\int\log\left(2\theta (\KE(2\theta)-\l)\right)d\hat \mu_{E_N}(\l)} \\
\int_{|\g_N|\ppq\e, |\hat \g_N|\ppq\e}
\exp\left\{ -\theta N \g_N {\hat \g_N - v\g_N\over \g_N+1}\right\}
\prod_{i=1}^N dP_i(g_i)
\end{multline} 

Now, under $\prod_{i=1}^N dP_i(g_i)$,
%$\left(\sqrt{N}\hat U_N,\sqrt{N}(\hat V_N - v\hat U_N)\right)
$\left(\sqrt{N}\g_N,\sqrt{N}\hat \g_N \right)
$ converges in law towards a centered two-dimensional
Gaussian variables $(\G_1, \G_2)$ as soon as their 
covariances converge. We investigate this convergence.

Hereafter, we shall write 
$g_i=(1+2\theta(v-\l_i))^{-{1\over 2}} \tilde g_i$
with standard independent Gaussian variables $\tilde g_i$.
Then,
$$
\E((\sqrt{N}\g_N)^2) = N\E\left[\left( {1\over N}\sum_{i=1}^N
\frac{ \tilde g_i^2-1}{1+2\theta v-2\t \l_i} +{1\over 2\theta}(H_{E_N}-\HE)(
\KE(2\theta))\right)^2\right]
$$
where we used that 
\begin{equation}
\label{nonews}
 2\theta=\HE(\KE(2\theta))=
\int \frac{1}{\KE(2\theta)-\l}d\mu_E(\l) \,,
\end{equation}
and (\ref{eqp}).
Equation (\ref{pesah}) implies
\begin{eqnarray*}
\lim_{N\ra\infty} \E((\sqrt{N}\g_N)^2)&=&
\lim_{N\ra\infty}
N\E\left[\left( {1\over N}\sum_{i=1}^N
\frac{\tilde g_i^2-1}{1+2\theta v-2\t \l_i} \right)^2\right]\\
&=&\lim_{N\ra\infty}
 {2\over N}\sum_{i=1}^N
\frac{1}{(1+2\theta v-2 \t \l_i)^2} \\
&=&{1\over 2\theta^2}\int \frac{1}{(\KE(2\theta)-\l)^2} d\mu_E(\l) := 
\frac{Z}{2\theta^2},
\end{eqnarray*}
where the above convergence holds since $\KE(2\theta)$ 
lies outside $[\l_{\rm min}, \l_{\rm max}]$ and therefore
outside the support of $\mu_E$. 

Similar computations give that under the same hypotheses,
$$
\lim_{N\ra\infty} \E((\sqrt{N}\hat \g_N)^2) =
 {1\over 2\theta^2}\int \frac{\l^2}{(\KE(2\theta)-\l)^2} d\mu_E(\l)\,$$
and that
$$
\lim_{N\ra\infty} \E(\sqrt{N}\hat \g_N\sqrt{N}\g_N)
= {1\over 2\theta^2}\int \frac{\l}{(\KE(2\theta)-\l)^2} d\mu_E(\l)\,.$$

Therefore, provided that the Gaussian integral is well defined, we find that
\begin{equation}
\label{pesahz}
I_N(\t,E_N)=e^{N\theta v-{N\over 2}
\int\log\left(2\theta(\KE(2\theta)-\l)\right)d\mun_{E_N}(\l)}
\int
e^{\{ -\theta x (y- v x)\}}
d\G(x, y)(1+o(1)),
\end{equation} 
with $\G$ a centered Gaussian measure on $\RR^2$ with covariance matrix
$$ R =\frac{1}{2\theta^2}\left[\begin{array}{ll}
\int \frac{1}{(\KE-\l)^2}d\mu_E(\l) &
 \int \frac{\l}{(\KE-\l)^2}d\mu_E(\l) \\
\int \frac{\l}{(\KE-\l)^2}d\mu_E(\l) &
\int \frac{\l^2}{(\KE-\l)^2}d\mu_E(\l) 
\end{array}
\right],$$
where we used the notation $\KE:= \KE(2\t).$

Following the ideas \cite{bolt}
as outlined in  appendix \ref{app}, we know that 
there is one step needed to justify this derivation, namely to check 
that the Gaussian integration in  (\ref{pesahz}) is non-degenerate.
If we set
$D:=4\theta^4 {\rm det} R$,
then, using the relation 
(\ref{nonews}),
one finds that $D=Z-4\theta^2$,  and that the Gaussian integral 
in (\ref{pesahz}) equals
$$\frac{\theta^2}{\pi \sqrt{D}}\int \exp\left(-\frac12
\sum_{i,j=1}^2  K_{i,j} x_i x_j\right) dx_1 dx_2\,,$$
where  the matrix $K$ equals $\theta \left[\begin{array}{ll}
-2v & 1 \\
1&0\\
\end{array}
\right] +R^{-1}$, that is
\begin{equation}
\label{defK}
K=
\frac{2\theta^2}{D}\left[\begin{array}{ll}
\int \frac{\l^2}{(\KE-\l)^2}d\mu_E(\l)-\frac{(\KE-\frac{1}{2\theta})D}
{\theta} &
-\int \frac{\l}{(\KE-\l)^2}d\mu_E(\l)+\frac{D}{2\theta} \\
-\int \frac{\l}{(\KE-\l)^2}d\mu_E(\l)+\frac{D}{2\theta} &
\int \frac{1}{(\KE-\l)^2}d\mu_E(\l) 
\end{array}
\right]
\end{equation}
Our task is to verify that $K$ is  positive definite.
It is enough to check that $K_{11}>0$
and ${\rm det}K>0$. Re-expressing $K_{11}$, one 
finds that
\begin{eqnarray*}
K_{11}&=&\frac{2\theta^2}{D}\left(
1-4\theta\KE+\KE^2 Z -\frac{1}{\theta}
(Z-4\theta^2)\left(\KE-\frac{1}{2\theta}\right)\right)\\
&=&\frac{2\theta^2}{D}\left( \left(\KE-\frac{1}{2\theta}\right)^2 Z
+\frac{Z}{4\theta^2} -1\right)\\
\end{eqnarray*}
But Schwarz's inequality applied to (\ref{nonews}) yields that
$Z>4\theta^2$ as soon as $\mu_E$ is not degenerate, 
implying that
$$K_{11}>
\left(\KE-\frac{1}{2\theta}\right)^2 Z\pgq 0\,,$$
as needed.
Turning to the evaluation of the determinant, note that
$$ {\rm  det} K
=\frac{4\theta^4}{D^2} Z\left(\frac{Z}{4\theta^2}-1\right)>0\,,$$
where the last inequality is again due to (\ref{nonews}).

$\bullet$ Let us finally consider the
case $\mu_E=\d_e$.
In this case, $\HE(x)=(x-e)^{-1}$ and 
$\KE(x)=x^{-1}+e$, $v=e$ (note also that $Z$
in Theorem \ref{precise1}.1 is equal to $4\theta^2$).
We can follow the previous proof but then
$$\lim_{N\ra \infty}
\E[(\sqrt{N}(\hat\g_N-v\g_N))^2]=0.$$
From here, we argue again using 
  appendix \ref{app} that
$$\lim_{N\ra\infty}
\E[1_{|\gamma_N|\le\e, |\hat\g_N-v\g_N|\le\e}\,\,
 e^{-\theta
 (1+\gamma_N)^{-1} \sqrt{N}\g_N \sqrt{N}(\hat\g_N-v\g_N)}
]=1$$
which completes the proof of Theorem \ref{precise1}. \hfill\xx

%%%%%%%%%%%%%%%%%%%%%%%%%%%%%%%%%%%%%%%%%%%%%%%%%%%%%%%%%%%%%%%%%%%%%%%%%%%%%
\section{Extension of the results to the complex plane}
\label{C}

In this section, we would like to extend the results of section
\ref{short} to the case where $\t$ is complex, that is to show
Theorem \ref{rg1complexe}.

As in the real case, we first would like to write that
\begin{equation}
\label{changevar}
I_N(\t, E_N) = \prod_{i=1}^N \sqrt{\zeta_i}
\int \exp\left\{\theta N{\sum_{i=1}^N\lambda_i\zeta_i g_i^2
\over \sum_{i=1}^N \zeta_i g_i^2} -{1\over 2}\sum_{i=1}^N
\zeta_i g_i^2\right\}\prod_{i=1}^N dg_i,
\end{equation}
with $\displaystyle \zeta_i = \frac{1}{1 +2\t v -2 \t \l_i},$ for 
$v$ such that $\Re(\zeta_i)>0,\quad \forall i \textrm{ with }1\le i\ppq N$.\\

This is a direct consequence of the following lemma
\begin{lem}
\label{CVg}
For any function $f: \CC^N \longrightarrow \CC$ which is invariant by $x \mapsto -x$, analytic outside $0$
and bounded on $\{z = x+iy\in \CC /  |y|< x \}^N$ and for any $(\z_1, \ldots, \z_N)$ such that $\Re(\z_i) > 0$ for any $i$ from $1$ to $N$, we have that
\begin{eqnarray*}
\mathcal J_N & := & \int f(g_1, \ldots, g_N) e^{-\frac 1 2  \sum_{i=1}^N g_i^2}  \prod_{i=1}^N dg_i \\
& = &  \prod_{i=1}^N \sqrt{\z_i} \int f( \sqrt{\z_1}g_1, \ldots,  \sqrt{\z_N}g_N) e^{-\frac 1 2  \sum_{i=1}^N \z_i g_i^2}  \prod_{i=1}^N dg_i,
\end{eqnarray*}  
with $\sqrt.$ is the principal branch of the square root in $\CC$.
\end{lem}

\textbf{Proof of Lemma \ref{CVg}:}

We denote by $r_j$ the modulus of $\z_j $ and 
$\a_j$ its phase $\left( \z_j = r_j e^{\a_j}\right)$.\\

As $f$ is bounded on $\RR^N$, dominated convergence gives that
$$\mathcal J_N = \lim_{R \ra \infty, \eps \ra 0} \int_{[-R,R]^N \setminus [-\eps,\eps]^N}
f(g_1, \ldots, g_N) e^{-\frac 1 2  \sum_{i=1}^N g_i^2}  \prod_{i=1}^N dg_i.$$
Thanks to invariance of $f$ by $x \mapsto -x$, we also have that
$$\mathcal J_N = \lim_{R \ra \infty, \eps \ra 0} 2^N \int_{[\eps,R]^N}
f(g_1, \ldots, g_N) e^{-\frac 1 2  \sum_{i=1}^N g_i^2}  \prod_{i=1}^N dg_i.$$

For each $j$ from $1$ to $N$ and $R \in \RR^+$, we define the following
segments in~$\CC$~:
$$ \mathcal C_{R, \eps}^{j} := \left\{ re^{i  \frac{\a_j}{2}}; \eps \ppq r \ppq R\right\},$$
and the following arc of circles
$$ \mathcal D_{ \eps}^{j} := \left\{ \eps e^{i  \a}; 0 \ppq \a \ppq  \frac{\a_j}{2} \right\} 
\textrm{ and } 
 \mathcal D_{R}^{j} := \left\{ R e^{i  \a}; 0 \ppq \a \ppq  \frac{\a_j}{2} \right\},$$
so that, for each $j$, $[\eps, R]$ run from $\eps$ to $R$
followed by $ \mathcal D_{R}^{j}$ run counterclockwise, followed by $ \mathcal C_{R, \eps}^{j}$
run from $Re^{i  \frac{\a_j}{2}}$ to $\eps e^{i  \frac{\a_j}{2}}$ followed by $ \mathcal D_{ \eps}^{j}$
run clockwise form a closed path.\\
Therefore, if we let
\[
f^{x_2, \ldots, x_N}_1 : \begin{array}[t]{rcl}
\CC & \rightarrow & \CC \\
x & \mapsto & f(x, x_2, \ldots, x_N),
\end{array}
\]
then for any $(x_2, \ldots, x_N) \in \CC^{N-1}$,
$x\mapsto f^{x_2, \ldots, x_N}_1(x)e^{-\frac 1 2 x^2}$ is analytic
inside the contour $[\e,R]\cup  \mathcal D_{R}^{j}\cup 
\mathcal C_{R, \eps}^{j}\cup \mathcal D_{\e}^{j}$,
so that Cauchy's theorem implies
\begin{multline*}
\int_{[\eps, R]} f^{g_2, \ldots, g_N}_1(g_1) e^{-\frac 1 2  g_1^2}  dg_1 = 
\int_{ \mathcal C_{R, \eps}^{1}}f^{g_2, \ldots, g_N}_1(g_1) e^{-\frac 1 2  g_1^2}  dg_1 \\
- \int_{ \mathcal D_{R}^{1}}f^{g_2, \ldots, g_N}_1(g_1) e^{-\frac 1 2  g_1^2}  dg_1 
+ \int_{ \mathcal D_{\eps}^{1}}f^{g_2, \ldots, g_N}_1(g_1) e^{-\frac 1 2  g_1^2}  dg_1 .
\end{multline*}
If we denote by 
$$ J_{N,R}^1 = \int_{[\eps, R]^{N-1}} e^{- \frac 1 2 \sum_{i=2}^N g_i^2}\int_{ \mathcal D_{R}^{1}}
f^{g_2, \ldots, g_N}_1(g_1) e^{-\frac 1 2  g_1^2}  dg_1 \ldots dg_N,$$
we have that
\begin{eqnarray*}
|J_{N,R}^1| & = &  \int_{[\eps, R]^{N-1}} \int_0^{\frac{\a_1}{2}} f(g_1, \ldots, g_N)
e^{- \frac 1 2 \sum_{i=2}^N g_i^2} R e^{- \frac 1 2 R^2 \cos(2 u_1)} du_1 dg_2 \ldots dg_N \\
& \ppq & \| f\|_\infty \sqrt{2 \pi}^N \frac{\a_1}{2} R e^{- \frac 1 2 R^2 \cos(\a_1)}.
\end{eqnarray*}
As $\cos(\a_1) >0$, we have that for any $\eps$, 
$\displaystyle  \lim_{R \ra \infty} |J_{N,R}^1| =0$.\\
In the same way, if we let
$$ L_{N,\eps}^1 = \int_{[\eps, R]^{N-1}} e^{- \frac 1 2 \sum_{i=2}^N g_i^2}\int_{ \mathcal D_{\eps}^{1}}
f^{g_2, \ldots, g_N}_1(g_1) e^{-\frac 1 2  g_1^2}  dg_1 \ldots dg_N,$$
then we have that
$$ |L_{N,\eps}^1| \ppq  \| f\|_\infty \sqrt{2 \pi}^N \eps \frac{\a_1}{2},$$
so that
$\displaystyle  \lim_{\eps \ra 0} |L_{N,\eps}^1| =0$.\\
By doing the same computation for each variable, we get that
\begin{multline*}
\lim_{R \ra \infty, \eps \ra 0}  \int_{[\eps,R]^N}
f(g_1, \ldots, g_N) e^{-\frac 1 2  \sum_{i=1}^N g_i^2}  \prod_{i=1}^N dg_i\\
=  \lim_{R \ra \infty, \eps \ra 0}  \int_{\prod_{i=1}^N \mathcal C_{R, \eps}^{1}}
f(g_1, \ldots, g_N) e^{-\frac 1 2  \sum_{i=1}^N g_i^2}  \prod_{i=1}^N dg_i.
\end{multline*}
The last step is to make the change of variable in $\RR$ which consist in 
letting $\tilde g_j = \sqrt{r_j} g_j$ to get the result announced in the 
lemma \ref{CVg} and therefore the formula (\ref{changevar}). \hfill \xx \\

We now go back to the \textbf{proof of Theorem \ref{rg1complexe}} and
proceed as in section \ref{short}. We let
 $$\g_N := {1\over N}\sum_{i=1}^N  \z_i g_i^2 - 1 \textrm{ and }
\hat \g_N := {1\over N}\sum_{i=1}^N\lambda_i \z_i g_i^2 -v(\t),$$
with $v(\t) = \RE(2\t)$,
which, for $|\t|$ small enough, is well defined
and such that $\Re \zeta_i >0$, by virtue of 
Property \ref{Rcomp} and Proposition \ref{Hcomp}.

Therefore, we find that
\begin{equation}
\label{f1}
I_N(\t, E_N)
= \prod_{i=1}^N\sqrt{\z_i} \, \, e^{N\theta v}
\int \exp\left\{ N\theta{\g_N(v\g_N-\hat \g_N)\over
1+\g_N}\right\} e^{ -{1\over 2}\sum_{i=1}^N
g_i^2}\prod_{i=1}^N dg_i,
\end{equation}
which is almost similar to what we got in (\ref{presque})
except that in the complex plane this is not so easy 
to ``localize'' the integral around $0$ as we did before.\\
 Our goal is now to show that 
$$ \lim_{N\ra\infty}
\int \exp\left\{ N\theta{\g_N(v(\theta)\g_N-\hat \g_N)\over
1+\g_N}\right\} e^{ -{1\over 2}\sum_{i=1}^N
g_i^2}\prod_{i=1}^N dg_i $$
exists and is not null.\\

Denote $\g_N=u^N_1+iu^N_2 - 1$
and $\hat \g_N=v^N_1+iv^N_2 - v(\t)$,
and let 
$$X^N + X_0:=(u^N_1,u^N_2, v^N_1,v^N_2) =
\left(\int \z_1(\l)x^2d\mun(x,\l),\dotsc,\int\z_4(\l)x^2
d\mun(x,\l)\right) $$
with $d\mun = \frac 1 N \sum_{i=1}^N \d_{\l_i,g_i},$
$$\z_1(\l)=\Re(  (1+2v(\theta) \theta-2\theta\l)^{-1}), \,\,
\z_2(\l)=\Im( (1+2v(\theta) \theta-2\theta\l)^{-1}),$$
$$
\z_3(\l)=\Re( \l (1+2v(\theta) \theta-2\theta\l)^{-1}), \,\,
\z_4(\l)=\Im( \l (1+2v(\theta) \theta-2\theta\l)^{-1})$$
$$\textrm{and } X_0=(1,0,\Re(v(\theta)),\Im(v(\theta))).$$
Then, 
we easily see as in \cite{BDG} (cf Lemma 4.1 therein)
 that the law of $X^N$ under $\prod_{i=1}^N \sqrt{2\pi}^{-1}
e^{-{1\over 2} g_i^2}dg_i$ satisfies a large deviation
principle on $\RR^4$ with rate function
$$\L^*(X)=\sup_{Y\in\RR^4\atop 1-2\langle \z(\l),Y 
\rangle\ge 0 \mu_E\mbox{ a.s. } }\left\{
\langle Y,X+X_0\rangle +{1\over 2}\int\log\left(1-2\langle \z(\l),Y 
\rangle \right)d\mu_E(\l)\right\},$$
with $\langle\, , \, \rangle$ the usual scalar product on $\RR^4$.\\
We denote
$$F(X^N):=\theta{\g_N(v\g_N-\hat \g_N)\over
1+\g_N}= F_1(X^N)+iF_2(X^N)$$
with $F_1$ and $F_2$ respectively the real and 
imaginary part of $F$. With these notations,  our problem boils down 
to show that
$
\E[e^{N F(X^N)}]$ converges towards a non-zero limit.
Following \cite{BA}, we know that it is  enough for us to
check that 
\begin{enumerate}
\item 
there is a vector  $X^*$ so that $F(X^*)=0$ and  
$$\lim_{M\ra\infty}\lim_{N\ra\infty}\left({1\over N}\log
\E[ e^{N  F_1(X^N)}]-{1\over N}\log
\E[ 1_{|X^N-X^*|\ppq \frac{M}{\sqrt{N}} }  e^{N  F_1(X^N)}]
\right)=0.$$
To prove this, the main part of the work will be
to show that\\
a) $X^*$ is the unique minimizer of
$\L^* - F_1$ (This indeed 
entails that the expectation can be localized in a small
ball around $X^*$), and then we will check that \\
b) $X^* $  is a not degenerate minimizer i.e 
the Hessian of $\L^*-F_1$ is positive
definite at $X^*$ (As shown in appendix \ref{app}, this will allow us to take this small ball of radius of order
$\sqrt{N}^{-1}$).
\item 
$X^*$ is also a critical point 
of $F_2$. This second point
allows to see that there is no
fast oscillations which reduces the
first order of the integral.
\end{enumerate}
Once these two points are checked,
it is not hard to see that

$$\E[ e^{N  F(X^N)}]= \E[e^{N D^2 F[X^*] (X^N-X,X^N-X)}](1+o(1))
=\mbox{det}(  D^2 (\Lambda^*-F)[X^*])^{-\frac 1 2}(1+o(1)).$$
This formula extends analytically
the result of Theorem \ref{rg1complexe}.  %Ader
In our case, $F$ depends linearly on $\theta$ 
and $X^*$ is the origin, from which it is 
easy to see that the convergence, if it holds
for some complex $\theta\neq 0$, will hold
in a neighborhood of the origin
since non degeneracy and
uniqueness of the minimizer 
questions will continuously depend on $\theta$. 
Moreover, it is not hard to see that
the convergence will actually 
hold uniformly in such a neighborhood 
of the origin (again because 
error terms will depend continuously on $\theta$).

$\bullet$ \textbf{Proof of the first point :}
To prove a),
let us notice that by our choice of $v(\theta)$ (see Proposition \ref{Hcomp}),
$\L^*$ is minimum 
at the origin and that the differential
of $F_1$ at the origin
is null. Hence, the origin is a critical
point of $F_1-\L^*$ (where this function is null)
and we shall now prove that it is the unique 
one when $|\theta|$ is small enough.\\
For that, we adopt the strategy used in  \cite{BDG}
and  consider
the joint deviations of the law of $(X^N,\mun)$.
A slight  generalization of Lemma 4.1
therein shows that it satisfies a
\ldp \, on $\RR^4 \ts \PR$ with good rate function
$$J(X,\mu)=I(\mu|\mu_E\otimes P)
+\tau \left(X+X_0-\int \z(\l) x^2 d\mu(\l,x)\right),$$
with $I(.|.)$ the usual relative entropy, $P$ a standard Gaussian measure
and
$$\tau(X)=\sup_{\a \in\Da_0}
\{\langle \a,X \rangle
\},
$$
where 
$\Da_0=\{\a\in\RR^4:1-2 \langle \a,\z(\l) \rangle \pgq 0
\quad \mu_E \textrm{ a.s. }
\}.$
From that and the contraction principle we have that
\begin{multline}
\label{po}
I(X):= \L^*(X)-F_1(X) \\
=\inf_{\mu\in\Pa(\RR)}\sup_{\a \in \Da_0}
\left\{I(\mu|\mu_E\otimes P)
+ \langle X+X_0-\int \z(\l) x^2 d\mu(\l,x),\a \rangle
-F_1(X) \right\} .
\end{multline}
If we set 
$$\mu^\a(dx,d\l)={1\over Z_\a}e^{-{1\over 2}(
1-2 \langle \z(\l),\a \rangle)x^2}
dx d\mu_E(\l)$$
then
$$I(\mu|\mu^\a)=I(\mu|\mu_E\otimes P)
- \langle \a,\int \z(\l)x^2d\mu(\l,x) \rangle 
-{1\over 2}\int\log(1-2 \langle \z(\l),\a \rangle)
d\mu_E(\l).$$
Thus,
\begin{multline*}
I(X)=\inf_{\mu\in\Pa(\RR)}\sup_{\a}\left\{
I(\mu|\mu^\a)+ \langle X+X_0,\a \rangle \right.\\
\left.+{1\over 2}\int\log(1-2 \langle \z(\l),\a \rangle)
d\mu_E(\l)-F_1(X)\right\}.
\end{multline*}
Observe that the supremum in $\L^*(X)$
is achieved at some $Y^X$
since $Y \mapsto -\int\log(1-2 \langle \z(\l),Y \rangle)
d\mu_E(\l)$ is lower semicontinuous
and $\{ Y\in \RR^4: 1-2\langle \z(\l),Y 
\rangle\ge 0\quad \mu_E\mbox{ a.s. }\}$ is compact
when $\mu_E$ is not a Dirac mass.
Indeed, from the definition of $v(\t)$,
we find that $\mu_E(\z_i(\l)>0)>0$ as
well as $\mu_E(\z_i(\l)<0)>0$ for $1\le i\le 4$ 
from which the compactness follows.
Moreover $Y^X$ satisfies 
\begin{equation}\label{eqY} (X+X_0)_i=\int
\frac{\z_i(\l)}{1-2<\z(\l),Y^X>}
d\mu_E(\l),\quad 1\le i\le 4.
\end{equation}
Consequently,
$$\L^*(X)-F_1(X)=I(X)\pgq \inf_{\mu\in\Pa(\RR)}\left\{
I(\mu|\mu^{Y^X}) +\L^*(X)-F_1(X)\right\}.$$
Since  $I(\mu|\mu^{Y^X})\pgq 0$, we deduce that
the infimum in $\mu$ is taken at $\mu=\mu^{Y^X}$. 
We  also check that  $\int \z(\l)x^2 d\mu^{Y^X}(\l,x)
=X + X_0$ due to (\ref{eqY}). Hence, going back to (\ref{po}), we find
that $I(X)=\Ia(\mu^{Y^X})$ with 
$$\Ia(\mu)=I(\mu|\mu_E\otimes P)
-F_1\left(\int \z(\l)x^2 d\mu(x,\l) - X_0\right).$$
We next show that $\Ia$ has a unique
minimizer for $\theta$
 small enough, and this  minimizer
satisfies
$\int \z(\l)x^2 d\mu(x,\l)= X_0$.
If the infimum is actually reached at a point $\mu^*$
such that $F_1$ is regular enough at the vicinity of 
$\int \z(\l)x^2 d\mu^*(x,\l) - X_0$ then this saddle point satisfy
the equation
\begin{equation}
\label{sp}
d\mu(x,\l)={1\over Z_\mu}
e^{DF_1(\int \zeta(\l)x^2 d\mu(x,\l)-X_0)[\zeta(\l)x^2]
-{1\over 2}x^2} dx d\mu_E(\l).
\end{equation}
Before going on the proof, let us justify that it is indeed
the case. Note first that 
as $\theta$ goes to zero,
 $v(\t)$ goes to $m=\int \l d\mu_E(\l)$
and $\Re[(1+2\theta v-2\theta\l)^{-1}]$ 
is bounded below by say $2^{-1}$.
Consequently,
 $\Re \g_N + 1 \pgq 2^{-1} {1\over N}\sum_{i=1}^N g_i^2$.
The rate function for the deviations of the
latest is $x-\log x-1$ which goes to infinity
as $x$ goes to zero as $\log x^{-1}$.
Therefore, for $\theta$ small enough,
$$\L^*(X)\pgq \log (2X_1)^{-1}$$
Since $F_1(X)$ is locally  bounded , we deduce that
the infimum has to be taken on $X_1\pgq \e$
for some fixed $\e>0$. In particular, $F_1$
is $\Ca^\infty$ on this set and equation (\ref{sp}) is
well defined.\\

We now want to use this saddlepoint equation to show 
uniqueness. Suppose that there are two minimizers $\mu$
and $\nu$ satisfying (\ref{sp}). Then
\begin{multline*}
\D:=\left|\int \z(\l)x^2 d\mu(x,\l)-
\int \z(\l)x^2 d\nu(x,\l)\right| \\ 
\ppq 4C|\theta| 
\sup_i  \int |\z_i(\l)|x^2
(d\mu(x,\l)+d\nu(x,l) ) \D , %Ader
\end{multline*}
as we have that $y\ra DF_1(y)[x]$
is Lipschitz, with Lipschitz norm 
of order $C|\theta| \|x\|$. 
We have now to show that for $\theta$ small 
enough, these covariances
are uniformly bounded. This can be done
using some arguments very similar to the ones 
we gave above to justify that the critical points
are such that $X_1\pgq \e$. We let it to the reader.
For $\theta$
small enough, we obtain a contraction so that
$\D =0$, which entails also $\mu=\nu$. It is easy to check that 
$\mu$ such that
$\int \zeta(\l)x^2 d\mu(x,\l) = X_0$ is always a solution to
(\ref{sp}),
and hence the unique one when $\t$ is small enough.
Observe now that  by (\ref{sp}),
this  minimizer is of the form $\mu^*=\mu^{\a^*}=\mu^{Y^{X^*}}$,
so that $X^*=\int \zeta(\l)x^2 d\mu^{\a^*}(x,\l) - X_0=0$
minimizes indeed $I$ and is actually its unique minimizer.

This concludes the proof of point a), which was the hard part of the work.\\

As we announced at the beginning and following \cite{BA}, we now have to show b),
that is to say to check that this minimizer
is non-degenerate. To see that, remark that
the second order
derivative of $F_1$ at the origin is
simply
\begin{equation}\label{Fo}
D^2F_1[0](U,V)=\Re(\theta(U(vU-V)))\ppq C|\theta|(|U|^2+|V|^2)
=C|\theta|\left(\sum_{i=1}^4 X_i^2\right) 
\end{equation}
with $U=X_1+iX_2,V=X_3+iX_4$.\\
\indent On the other side, 
observe that, as $d(\mun_{E_N}, \mu_E) = o(\sqrt N^{-1})$, the covariance 
 matrix of \linebreak
$\sqrt N (u_1^N, (\Im(\theta))^{-1} u^N_2,
v^N_1, (\Im(\theta))^{-1} v^N_2)$ converges as $N$
goes to infinity towards a $4\ts4$ matrix $K(\theta)$
which is positive definite. 
Now, remark that $v(\t)=\RE(2\theta)$ implies that
$\Re(\theta)(\Im(\theta))^{-1} \Im(v(\theta))$ converges
as $|\theta|$ goes to zero, from which we argue
that $K(0)$ is positive definite and bounded. By continuity in $\theta$
of $K(\theta)$ we deduce that $ K(\theta)\le C I$
for some $C>0$ and $\t$ small enough.
and the limiting 
covariances $\sqrt{N}(u^N_1,  u^N_2,
v^N_1,  v^N_2)$(which are also given by the second order derivatives
of $\L^*$)
 converges towards  a matrix $K^\prime(\theta)$ such that
$$ D^2\L^*[0](X,X)=
\langle X,K^\prime(\theta)^{-1}X \rangle \ge C^{-1}(X_1^2+X_3^2+(\Im(\theta))^{-2}X_2^2
+(\Im(\theta))^{-2}X_4^2)$$
and hence, this together with \eqref{Fo} gives that, 
for $|\t|$ small enough, $\frac 1 2 D^2\L^*[0] - D^2F_1[0] \pgq 0. $

$\bullet$ \textbf{Proof of the second point :}
To get Theorem \ref{rg1complexe},
the last step is now to establish the second point, namely to check that $0$ is also a critical point for
$F_2$, which is  straightforward computation 
since $F$ behaves in the neighborhood
of the origin as a sum of monomials
of degree $2$ in $X$.
\hfill \xx

%%%%%%%%%%%%%%%%%%%%%%%%%%%%%%%%%%%%%%%%%%%%%%%%%%%%%%%%%%%%%%%%%%%%%%%%%%%%%%%%%%%%%%%%%%%%%%%%%%%%%%%%%%%%%%%%%%%%%%%%%%%%%%%%%%%%%%%%%%%%%%%%%%%%%
\section{Full asymptotics in the real rank one case}
\label{rang1}
The goal of this section is to establish the convergence and to 
find an explicit expression for 
$\displaystyle \IE(\t) := \lim_{N \ra \infty} 
\frac 1 N \log I_N(\t, E_N)$ 
as far as $E_N$ satisfies Hypothesis \ref{hypsurEN}
but $\t$ do not necessarily satisfy the hypotheses of Theorem \ref{sh}.
This corresponds to show
Theorem \ref{expreel} (we again restrict to the case $\b = 1$ to
avoid heavy notations).\\
We recall that 
$$ I_N(\t, E_N) = \EE\left[ \exp\left(N\t \frac{\sum_{i=1}^N \l_i g_i^2}{\sum_{i=1}^N g_i^2}\right)\right],$$
therefore one main step of the proof will be to get a large deviation principle for $\displaystyle z_N:=\frac{\sum_{i=1}^N \l_i g_i^2}{\sum_{i=1}^N g_i^2}$.

\subsection{Large deviation
bounds for $z_N$}

We denote by $\displaystyle u_N := \frac 1 N \sum_{i=1}^N g_i^2$
and  $\displaystyle v_N := \frac 1 N \sum_{i=1}^N \l_i g_i^2.$
We intend to get the following result
\begin{prop}
\label{gdz}
If the empirical measure $\mun_{E_N} = \frac 1 N \sum_{i=1}^N \d_{\l_i}$ satisfies
Hypothesis  \ref{hypsurEN}, 
the law $\hat \pi_N$ of 
$\left( u_N^{-1} v_N \right)$ 
under the standard $N$-di\-men\-sio\-nal Gaussian measure satisfies a
large deviation principle in the scale $N$ with good rate function
$$T (\a)=\left\lbc
\begin{array}{ll}
\frac 1 2 h_\a(\KE(\QE(\a))) & \mbox{ if }
\a \in [\a_{\rm min}, \a_{\rm max}],\\
\frac 1 2  h_\a^{\rm max} & 
\mbox{ if }\a \in ]\a_{\rm max}, \lM[,  \\
\frac 1 2  h_\a^{\rm min} &
\mbox{ if }\a \in ]\lm, \a_{\rm min}[ , \\
 +\infty
\mbox{ if } \a \notin ]\lm, \lM[\\
\end{array}
\right. $$
with $$ \a_{\rm max} = \lM - \frac 1 \HM \quad \textrm{ and }\quad  \a_{\rm min} = \lm - \frac 1 \Hm,$$
where  we recall that $\HM = \lim_{z \downarrow \lM} \int \frac{1}{z-\l} d\mu_E(\l)$ and 
$\Hm = \lim_{z \uparrow \lm} \int \frac{1}{z-\l} d\mu_E(\l)$; \linebreak 
we denote also,
for $\k \in [\lm, \lM]^c$,
$$ h_\a(\k) = \int \log\left( \frac{\k - \l}{\k-\a}\right) d\mu_E(\l),$$ 
$h_\a^{\rm min} = \lim_{\k \uparrow \lm} h_\a(\k)$ and 
$h_\a^{\rm max} = \lim_{\k \downarrow \lM} h_\a(\k)$.
Finally, the functions $\KE$ and $\QE$ were defined respectively in Definition 
\ref{KRreel} and Property \ref{Rreel}. \\
Note that $\HM$ and $\Hm$ can be infinite (respectively $+\infty$
and $-\infty$); in this case,  we adopt the convention that $\frac{1}{\infty} =0$.
\end{prop}

The \textbf{proof of Proposition \ref{gdz}} decomposes mainly in four steps,
expressed in the following four lemmata :
\begin{lem}
\label{lem1}
  For any $\a \in [\lm, \lM]$,
\begin{align*}
 \lim_{\e\ra 0}\varliminf_{N\ra\infty} &\frac 1N \log
\hat\pi_N\left( \bigl|{v_N}  -\a {u_N}\bigr|<\sqrt{\e}\right)  \\
 & \ppq \lim_{\e\ra 0}\varliminf_{N\ra\infty}\frac 1N \log
\hat\pi_N\left( \left|\frac{v_N}{u_N} -\a\right|<\e\right)   \\
& \ppq \lim_{\e\ra 0}\varlimsup_{N\ra\infty}\frac 1N \log
\hat\pi_N\left(  \left|\frac{v_N}{u_N} -\a \right|<\e\right)  \\
& \qquad \qquad \ppq \lim_{\e\ra 0}\varlimsup_{N\ra\infty}\frac 1N \log
\hat\pi_N\left( \bigl|{v_N} -\a {u_N}\bigr|<\sqrt{\e}\right) \label{encadr} 
\end{align*}
\end{lem}

\begin{lem}
\label{lem2}
We denote by $v_N(\g) :=N^{-1}\sum_{i=1}^N \g_i g_i^2$ and we assume that
the $\g_i$'s are such that
\begin{enumerate}
\item $\g_{\rm max}^N:=\max_{1\le i\le N} \g_i$
(resp. $\g_{\rm min}^N=\min_{1\le i\le N} \g_i$)
converges towards $\g_{\rm max}<\infty$ (resp. $\g_{\rm
min}>-\infty$).
\item The empirical measure $N^{-1}\sum_{i=1}^N \d_{\g_i}$ converges to 
a compactly supported measure $\mu$; we denote by $\g^+$ and $\g^-$
the edges of the support of $\mu$.
\end{enumerate}
Then, the law of $v_N(\g)$ satisfies a
large deviation principle in the
scale $N$ with rate function
$$J_{\mu,\g_{\rm min}, \g_{\rm max}} (x)=\left\lbc
\begin{array}{ll}
L(x)& \mbox{ if }
x\in [x_1, x_2]\\
L(x_1) +\frac{1}{2 \g_{\rm min}}(x-x_1)& 
\mbox{ if } x<x_1\\
L(x_2) +\frac{1}{2 \g_{\rm max}}(x-x_2) &
\mbox{ if } x>x_2\\
\end{array}
\right.$$
with
$$L(x)=\sup\left\{ ux+\frac 12 \int\log(1-2\l u) d\mu(\l)\right\}$$
where the supremum is taken over $u$ such that
$1-2\l u>0$ for every 
$\l \in [\g_{\rm min},\g_{\rm max}]$,
\[
x_1 = \left\{ \begin{array}{ll}
 \g_{\rm min}( \g_{\rm min} \Hm^\g - 1), & \textrm{ if } \g_{\rm min} < 0\\
-\infty & \textrm{ otherwise, }
\end{array}\right.
\]
whereas
\[
x_2 = \left\{ \begin{array}{ll}
 \g_{\rm max}( \g_{\rm max} \HM^\g - 1), & \textrm{ if } \g_{\rm max} > 0\\
\infty & \textrm{ otherwise,}
\end{array}\right.
\]
with the obvious notations $\HM^\g = \lim_{z\downarrow \g_{\rm max}} H_\mu(z)$
and $\Hm^\g = \lim_{z\uparrow \g_{\rm min}} H_\mu(z)$.

\end{lem}

\begin{lem}
\label{lem3}
If we denote $\g_i^\a := \l_i - \a$, $\mu^\a$ the weak limit of the
empirical measure $\frac 1 N  \sum_{i=1}^N \d_{\g_i^\a}$ (note that
$\mu^\a$ is just $\tau_{-\a} \sharp \mu$, where $\tau_{-\a} $ is the shift
given by $\tau_{-\a}(x) = x-\a$), $\g_{\rm max}^\a$ and $\g_{\rm min}^\a$
are respectively the limits of ${\rm max} \g_i^\a$ and ${\rm min} \g_i^\a$,
then
$$ J_{\mu^\a, \g_{\rm max}^\a, \g_{\rm min}^\a}(0) = T(\a),$$
with $T $ as defined in Proposition \ref{gdz}.
\end{lem}

\begin{lem}
\label{lem4}
$T$ is a good rate function.
\end{lem}

Then,  Proposition \ref{gdz} follows easily from these lemmata.
 Indeed,by definition of $u_N$ and $v_N$,
 we have that, for all $\e>0$ and $N$
large enough  $z_N \in [\lm-\e, \lM +\e]$ 
so that,
$$ \limsup_{N\ra\infty}\frac 1N \log 
\hat\pi_N\left( z_N \in [\lm -\e,\lM +\e]^c\right) = -\infty. $$ 
Thus,
from Theorem 4.1.11 in \cite{DZ}, it is enough to consider small balls
ie to show that, for any $\a\in [\lm,\lM]$,
$$ \limsup_{\e\ra 0}\limsup_{N\ra\infty}\frac 1N \log 
\hat\pi_N\left( \bigl|z_N -\a \bigr|\le \e\right) \ppq -T(\a),$$
and 
$$ \liminf_{\e\ra 0}\liminf_{N\ra\infty}\frac 1N \log 
\hat\pi_N\left( \bigl|z_N -\a \bigr|<\e\right) \pgq -T(\a).$$
Now, if $\g_i^\a = \l_i - \a$ and the $\l_i$'s satisfy Hypothesis 
\ref{hypsurEN}, $v_N(\g^\a):= \frac 1 N \sum(\l_i-\a) g_i^2 = v_N-\a u_N$
satisfy the hypotheses $(1)$ and $(2)$ of Lemma \ref{lem2}.
Therefore it satisfies a large deviation principle with rate function
$ J_{\mu^\a, \g_{\rm max}^\a, \g_{\rm min}^\a}$.
In particular this gives that in  Lemma \ref{lem1},
the rightmost and leftmost members coincide, so that
\begin{multline*}
\lim_{\e\ra 0}\liminf_{N\ra\infty}\frac 1N \log
\hat\pi_N\left( \bigl|{v_N} -\a {u_N}\bigr|<\sqrt{\e}\right) \\
= \lim_{\e\ra 0}\limsup_{N\ra\infty}\frac 1N \log
\hat\pi_N\left( \bigl|{v_N} -\a {u_N}\bigr|<\sqrt{\e}\right) 
= -J_{\mu^\a,\g_{\rm min}^\a, \g_{\rm max}^\a}(0) = -T(\a)
\end{multline*}
where the last equality comes from Lemma \ref{lem3}.\\
The study of the function $T$, that will give Lemma \ref{lem4},
allows to conclude the proof.

\subsection{Proofs of the lemmata}

\textbf{Proof of Lemma \ref{lem1}:}\\
For any $\a \in\RR$ and $\e >0$, we have 
\begin{multline*}
\hat\pi_N\left( |{v_N} -\a {u_N}|<\sqrt{\e}\right)-
\hat\pi_N\left( |u_N|\pgq \sqrt{\e}^{-1}\right)
\ppq \hat\pi_N\left( \left|\frac{v_N}{u_N} -\a\right|<\e\right) \\
\ppq
\hat\pi_N\left( |{v_N} -\a{u_N}|<\sqrt{\e}\right)+
\hat\pi_N\left( |u_N|\pgq \sqrt{\e}^{-1}\right).
\end{multline*}
Now, by Chebychev's inequality,
$$\hat\pi_N\left( |u_N|\pgq \sqrt{\e}^{-1}\right)
\ppq e^{-\frac{1}{4\sqrt{\e} }N} \hat\pi_N\left(e^{\frac{1}{4}
u_N}\right)
\ppq 2^N e^{-\frac{1}{4\sqrt{\e} }N},$$
so that
$$ \lim_{\e\ra 0}\limsup_{N\ra\infty}\frac 1N \log 
\hat\pi_N\left( |u_N|\pgq \sqrt{\e}^{-1}\right) = -\infty,
$$
what gives immediately Lemma \ref{lem1}.\\

Lemma \ref{lem2} is proved  in \cite{BGR},  Theorem 1;
we omit it here.

\textbf{Proof of Lemma \ref{lem3}:}\\   
Our goal is to identify $T(\a) = J_{\mu^\a,\g_{\rm min}^\a, \g_{\rm max}^\a}(0)$. As we said above, it is enough to restrict to 
$\a \in [\l_{\rm min}, \l_{\rm max}]$.\\
We have of course $\g_{\rm min}^\a = \l_{\rm min}-\a$ and 
$\g_{\rm max}^\a = \l_{\rm max}-\a$ and it is easy to check that
$$\HM^\a := \lim_{z \downarrow \lM -\a} \int \frac{1}{z-\l} d\mu^\a(\l) = \HM$$
(and respectively for $\Hm$).\\
Therefore, if we denote by $x_1^\a$ and $x_2^\a$ the bounds corresponding
to $\mu^\a$, we have that :
$$ x_1^\a = (\lm -\a) ((\lm -\a)\Hm -1)$$ 
(as the inequality $\g_{\rm min}^\a = \lm-\a <0$ is always satisfied for the
$\a$'s we are interested in) and similarly
$ x_2^\a = (\lM -\a) ((\lM -\a)\HM -1).$
We now have to determine the sign of $x_1^\a$ and $x_2^\a$ with respect to
$\a$. It is easy to check that
\begin{itemize}
\item $x_1^\a \ppq 0$ and $x_2^\a \pgq 0$ if $\displaystyle \a \in \left[\a_{\rm min}:= \lm - \frac  1 \Hm, \a_{\rm max}:= \lM - \frac  1 \HM \right]$ 
\item $x_1^\a \ppq 0$ and $x_2^\a \ppq 0$ if $\displaystyle \a \in [\a_{\rm max}, \l_{\rm max}]$
\item  $x_1^\a \pgq 0$ and $x_2^\a \pgq 0$ if $\displaystyle \a \in [\a_{\rm min}, \l_{\rm min}]$
\end{itemize}
Therefore, we deduce
$$J_{\mu^\a,\g_{\rm min}^\a, \g_{\rm max}^\a}(0)=\left\lbc
\begin{array}{ll}
L^\a(0) & \mbox{ if }
\a \in [\a_{\rm min}, \a_{\rm max}]\\
L^\a(x_2^\a) -\frac{1}{2} \HM(\a_{\rm max}-\a) &
\mbox{ if } \a_{\rm max} \ppq \a \ppq \lM\\
L^\a(x_1^\a) -\frac{1}{2} \Hm(\a_{\rm min}-\a) &
\mbox{ if } \lm \ppq \a \ppq \a_{\rm min},\\
\end{array}
\right.$$
where we recall that 
$$ L^\a(x) = \sup \left\{ ux + \frac 1 2 \int \log(1+2\a u - 2\l u) d\mu(\l)\right\},$$
with the supremum on $u$ such that $1+2\a u - 2\l u>0$ 
for all $\l \in [\l_{\rm min}, \l_{\rm max}]$.\\

We now get interested in the expression of $L^\a $ on $[x_1^\a, x_2^\a]$.\\
Obviously, the supremum is not reached at $u=0$.\\
For $u \neq 0$, we denote $\k := \a + \frac{1}{2u}$,
then we have that $1+2\a u - 2\l u = \frac{\k-\l}{\k - \a}$.
Moreover, if for all $\l \in [\l_{\rm min}, \l_{\rm max}]$,
$1+2\a u - 2\l u>0$ then ($\k > \lM$ and $u>0$) or ($\k < \lm$ and $u<0$) 
and conversely, so that

$$ L^\a(x)  =  \frac 1 2 \sup_{\k \in [\l_{\rm min}, \l_{\rm max}]^c}
\left\{  \frac{x}{\k - \a} + h_\a(\k)\right\},$$
with the notations of Proposition \ref{gdz}.\\

$\bullet$ If $\a \in \Ipp := [\a_{\rm min}, \a_{\rm max}]$,
$$ J_{\mu_\a,\g_{\rm min}^\a, \g_{\rm max}^\a}(0) = L^\a(0)
= \frac 1 2 \sup_{\k \in [\l_{\rm min}, \l_{\rm max}]^c} h_\a(\k).$$
We now want to check that in this case, the supremum of $h_\a$
is reached at $\k_0 = \KE(\QE(\a))$.\\

The first point is to show that in this case, there is a unique 
$\k_0$ where $h_\a^\prime$ cancels. Indeed :
$$ h^\prime_\a(\k_0)=0 \Longleftrightarrow \HE(\k_0) = \frac{1}{\k_0 - \a}
\Longleftrightarrow \k_0 = \KE(\QE(\a))$$
We now check that the maximum of $h_\a$ is reached at 
$\k_0$; 
\begin{itemize}
\item if $\k_0 > \lM$, $h_\a$ is decreasing from $0$ to $h^\a_{min}$ 
on $]-\infty, \lm[$, it is increasing from $h^\a_{max}$ to $h_\a(\k_0)$
on $]\lM,\k_0 ]$ and then decreasing from $h_\a(\k_0)$ to $0$ 
on $]\k_0, +\infty ]$,
\item if $\k_0 < \lm$, $h_\a$ is increasing from $0$ to $h_\a(\k_0)$
on $]-\infty, \k_0]$ then decreasing from $h_\a(\k_0)$ to $h^\a_{min}$ 
on $]\k_0, \lm[$, it is increasing from $h^\a_{max}$ to $0$ 
on $]\lM, +\infty [$.
\end{itemize}
We treat in details the proof of the first point, when $\k_0 > \lM$, 
the other one being very similar.
We recall from Property \ref{Rreel} that $\Ipp$ is the image of $\RE$.\\
If $\k_0 > \lM$, $h_\a^\prime$ does not cancel on $]-\infty, \lm[$.
It is negative since, when $\a \in \Ipp$, $ \lm-\frac 1 \Hm$
and so  $\lim_{\k \ra \lm} h_\a^\prime(\k)  < 0$.
On the other side, we want to find the sign of $h_\a^\prime$
on $]\lM, +\infty [$ knowing that it cancels at $\k_0$. 
As above, we show that 
$ \lim_{\k \ra \lM} h_\a^\prime(\k) > 0$ and we deduce from that 
and the continuity of $h_\a^\prime$,
that  it is positive till $\k_0$. 
Furthermore, $h_\a$ is also twice differentiable at $\k_0$ and
\begin{eqnarray*}
h_\a^{\prime\prime}(\k_0) & = & - \int \frac{1}{(\k_0 -\l)^2} d\mu_E(\l) 
+ \left(\frac{1}{\k_0-\a}\right)^2 \\
& < & - \left(\int \frac{1}{\k_0 -\l} d\mu_E(\l) \right)^2 
+(\HE(\k_0))^2 <0, 
\end{eqnarray*}
where we used Cauchy-Schwarz inequality and the definition of $\k_0$.
Therefore $h_\a^\prime$
is negative for $\k > \k_0$
and the fact that $ \lim_{\k \ra +\infty} h_\a(\k) =0$
concludes the proof of the first point.\\

Finally, we got that if $\a \in  [\a_{\rm min}, \a_{\rm max}]$,
$$ J_{\mu^\a,\g_{\rm min}^\a, \g_{\rm max}^\a}(0) = \frac 1 2 h_\a (\KE(\QE(\a)))$$

$\bullet$ If $\a > \a_{\rm max}$, our starting point is
$$ J_{\mu^\a,\g_{\rm min}^\a, \g_{\rm max}^\a}(0) = \frac 1 2
\sup_{\k \in [\lm, \lM]^c} \left\{ x_2^\a \left( \frac{1}{\k-\a}
- \frac{1}{\lM-\a} \right) + h_\a(\k)\right\}$$
Using arguments as above, we show that the function
$$ g_\a(\k) = \frac{x_2^\a}{\k-\a}  + h_\a(\k)$$
on $ [\lm, \lM]^c$ takes its supremum as $\k$ goes to $\lM$
by showing that its derivative is negative on $[\lm, \lM]^c$.
Hence, $ J_{\mu^\a,\g_{\rm min}^\a, \g_{\rm max}^\a}(0) = 
\frac 1 2 h_{\rm max}^\a$. 

$\bullet$ The case $\a< \a_{\rm min} $ is treated similarly,
which concludes the proof of Lemma \ref{lem3}.\\

The \textbf{proof of Lemma \ref{lem4}} is easy : $T$
is in fact continuous on $]\lm, \lM[$. 
Indeed, it is continuous on each interval $]\lm, \a_{\rm min}[$,
$]\a_{\rm min}, \a_{\rm max}[$ and $]\lM, \a_{\rm max}[$  so that it is enough to check
that $\KE(\QE(\a)) \xrightarrow[\a \ra \a_{\rm max}]{} \lM$ (see Property
\ref{Rreel}) so that $T(\a) \xrightarrow[\a \ra \a_{\rm max}]{} \frac 1 2
h^\a_{\rm max}$; and similarly at $\a_{\rm min}$.\\

\subsection{Proof of Theorem \ref{expreel}}

By Varadhan's lemma, we have 
\begin{lem}
\label{lem5}
For any $\t \in \RR$, if $T$ is the function introduced
in Proposition \ref{gdz}, we have
$$ \lim_{N \ra \infty} \frac 1 N \log I_N(\t, E_N)
= \sup_\a \{\t \a - T(\a)\}.$$
\end{lem}

Lemma \ref{lem5} therefore gives the existence of the limit,
the last step 
to conclude the proof of Theorem \ref{expreel} is to check that it coincides with the function $\IE$
introduced in Theorem \ref{expreel}. \\

We denote by 
$$ G(\t):= \sup_{\a \in \Ipp} \left[ \t \a - \frac 1 2 h_\a(\KE(\QE(\a)))\right],$$
$$ G_1(\t):= \sup_{\a \in I_1} \left[ \t \a - \frac 1 2 h^\a_{\rm max} \right], \quad
G_2(\t):= \sup_{\a \in I_2} \left[ \t \a - \frac 1 2 h^\a_{\rm min} \right],$$
where we recall that $\Ipp = [\a_{\rm min}, \a_{\rm max}]$
and we denote by $I_1 = ]\a_{\rm max}, \lM]$ and $I_2 = [\lm, \a_{\rm min}[$.\\ 

The main part of the work for this last step will rely on proving
\begin{lem}
\label{exprG}
With the notations introduced above, we have\footnote{$\sharp = - \infty$ if $\Hm = - \infty$ and otherwise these expressions 
are well defined in virtue of the fact that
$\displaystyle \int_0^1 \frac 1 \l d\mu(\l) < +\infty \Rightarrow  -\int_0^1 \log \l d\mu(\l) < +\infty$,\\
$* = - \infty$ if $\HM = + \infty$ and otherwise these expressions 
are well defined for the same reason.}
\[
 G(\t) = \left\{
\begin{array}[c]{ll}
\frac 1 2 \int_0^{2\t} \RE(u)du, & \textrm{ if } 2\t \in \Ip\cup \{0\} = ]\Hm, \HM[ \\
\t \a_{\rm min} - \frac 1 2 \int \log (\Hm(\lm-\l))d\mu_E(\l)^\sharp, & \! \!\textrm{ if } 2\t  \!\ppq  \!\Hm \\
\t \a_{\rm max} - \frac 1 2 \int \log (\HM(\lM-\l))d\mu_E(\l)^*, &   \! \! \!\textrm{ if } 2\t  \! \pgq  \!\HM ,
\end{array} 
\right.
\]

\[
 G_1(\t) = \left\{
\begin{array}[c]{ll}
\t\left( \lM - \frac{1}{2 \t} \right) - \frac 1 2 \int \log (2\t (\lM-\l))d\mu_E(\l)^*, &   \!\!\textrm{ if } 2\t  \! >  \!\HM \\
\t \a_{\rm max} - \frac 1 2 \int \log (\HM(\lM-\l))d\mu_E(\l)^*, &  \! \!\textrm{ if } 2\t  \!<  \!\HM ,
\end{array} 
\right.
\]

\[
 G_2(\t) = \left\{
\begin{array}[c]{ll}
\t\left( \lm - \frac{1}{2 \t} \right) - \frac 1 2 \int \log (2\t (\lm-\l))d\mu_E(\l)^\sharp, &   \!\! \!\textrm{ if } 2\t  \! <  \!\Hm \\
\t\left( \lm - \frac 1 \Hm \right) - \frac 1 2 \int \log (\Hm(\lm-\l))d\mu_E(\l)^\sharp, &   \!\! \!\textrm{ if } 2\t  \!>  \!\Hm .
\end{array} 
\right.
\]
\end{lem}

\noindent \textbf{Proof of Lemma \ref{exprG} :}\\
$\bullet$ We first study $G$. \\
This is finding the supremum of
$j_\t(\a) := \t \a - \frac 1 2 h_\a(\KE(\QE(\a))) $ on $\Ipp$.
From Definition \ref{KRreel} and Property \ref{Rreel}, we have that 
$j_\t$ is differentiable  on $\Ipp$ and an easy computation gives
$$ j_\t^\prime(\a) = \frac 1 2 (2 \t - \QE(\a)).$$ 
\begin{itemize}
\item If $2 \t \in \Ip$, $j_\t$ is maximized at
$\a_0 = \RE(2 \t)$ 
and so,  if $2\t \in ]\Hm, \HM[ \setminus \{0\}$,
\begin{eqnarray*}
 G(\t) &=& \frac 1 2 \left(2\t \RE(2\t) - \log(2\t)
- \int \log (\KE(2\t) - \l)d\mu_E(\l)
\right)\\
&=&= \frac 1 2 \int_0^{2\t} \RE(u) du. \\
\end{eqnarray*}

\item If $\Hm > - \infty$ and $2 \t < \Hm$, the 
equation $  j_\t^\prime(\a_0) = 0$ has no solution
and actually  $  j_\t^\prime$ is negativeso that  the supremum is reached at the left boundary 
$\displaystyle \a_{\rm min}$ of $\Ipp$
and is equal to 
$$ \t \a_{\rm min} - \frac 1 2 \int \log (\Hm(\lm-\l))d\mu_E(\l).$$

\item If $\HM < + \infty$, a similar treatment in the case $2\t > \HM$ concludes the proof for $G$.
\end{itemize}

$\bullet$ The formulas for $G_1$ and $G_2$
are derived similarly.

By virtue of Lemmata \ref{lem5} and \ref{exprG},
 to finish the proof of Theorem \ref{expreel}, 
we have now 
\begin{enumerate}
\item to compare $G_{|\Ip}$,  $G_{1|\Ip}$ and  $G_{2|\Ip}$ to get  $I_{\mu_E |\Ip}$.\\
Since 
$\displaystyle \lim_{\a \uparrow \HM} j_\t(\a) = G_1(\t)$ and
$\displaystyle \lim_{\a \downarrow \Hm} j_\t(\a) = G_2(\t)$
whereas $G(\t) = \sup_{\a \in \Ip} [j_\t(\a)]$, we get that
 $I_{\mu_E |\Ip} = G_{|\Ip}$.\\
\item if $\HM < +\infty$, to compare $G_{|\{ 2\t > \HM\}}$, $G_{1|\{ 2\t > \HM\}}$ 
and $G_{2|\{ 2\t > \Hm\}}$ to get $I_{\mu_E ||\{ 2\t > \HM\}}$.\\
By studying the function 
$\displaystyle x \mapsto - \frac{\t}{x} - \frac 1 2 \log x $, 
which reaches its maximum at $\t$, we can easily deduce that
 $G_{|\{ 2\t > \HM\}} < G_{1|\{ 2\t > \HM\}}$.\\
Moreover $G_{1|\{ 2\t > \HM\}}$ and  $G_{2|\{ 2\t > \HM\}}$ 
are the limits of $j_\t$ respectively at 
$\a_{\rm max}$ and  $\a_{\rm min}$
and we know that in the case $2\t > \HM$, $j_\t$ is
increasing. This gives $G_{2|\{ 2\t > \HM\}} < G_{1|\{ 2\t > \HM\}}$.\\
In this case we conclude that the maximum is given by 
$G_{1|\{ 2\t > \HM\}}$.
\item Arguing similarly, we can see that in the case where 
$ 2\t < \Hm$ the maximum is given by  $G_{2|\{ 2\t < \Hm\}}$. 
\end{enumerate} 

To conclude the proof of Theorem \ref{expreel},
we use the continuity of $\IE$ with respect to $\t$ given by
the first point of Lemma \ref{cont} 
to specify its value at $\lm$, $\a_{\rm min}$,
$\a_{\rm max}$ and $ \lM $. \hfill \xx

%%%%%%%%%%%%%%%%%%%%%%%%%%%%%%%%%%%%%%%%%%%%%%%%%%%%%%%%%%%%%%%%%%%%%%%%%%%%%%%%%%%%%%%%%%%%%%%%%%%%%%%%%%%%%%%%%%%%%%%%%%%%%%%%

\section{Asymptotic independence and free convolution}
\label{add}

In this section,
we want to prove Theorem \ref{Fourier},
that is to say  concentration and decorrelation properties for
the spherical integrals. 

We recall first that as an immediate Corollary 
of Theorem \ref{Fourier}, we
get that 
\begin{cor}\label{additive}
For $\theta$ sufficiently small
$$R_{\mu_B\boxplus \mu_A}(\theta)=R_{\mu_A}(\theta)
+R_{\mu_B}(\theta),$$
where $\boxplus$ denotes the free convolution of measures.
\end{cor}
\prf \
In fact, being given $\mu_A$, $\mu_B$,
we take $\l_1(A)$ (resp. $\l_1(B)$) to be the 
lower edge of the support of $\mu_A$ (resp. $\mu_B$)
and then set for $i\ge 2$
\begin{align*}
\l_i(A) &=\inf\left\{ x\ge \l_{i-1}(A): \mu_A([\l_1(A),x])\ge
{i\over N} \right\}, \\
\l_i(B) & =\inf\left\{ x\ge \l_{i-1}(A): \mu_B([\l_1(B),x])\ge
{i\over N}\right\}.
\end{align*}
It is easily seen that with this choice,
$A_N=\diag(\l_i(A))$ and $B_N=\diag(\l_i(B))$
satisfy Hypothesis \ref{hypsurEN}.
Since $\mu_A$ and $\mu_B$ are compactly supported,
$A_N$ and $B_N$
have uniformly bounded spectral radius and
so does $A_N+UB_NU^*$.
Hence, for $\theta$ small enough,
$A_N$, $B_N$ and $A_N+UB_NU^*$
satisfy the hypotheses of Theorem \ref{sh}
(recall that $A_N$ and $UB_NU^*$
are asymptotically free (c.f Theorem 5.2 in \cite{BP})
so that $\mun_{A_N+UB_NU^*}$ converges towards $\mu_B\boxplus \mu_A$).
Moreover, we can check that
$d(\mun_{A_N}, \mu_A) \ppq 2\|A_N\|_\infty \,N^{-1}$
and similarly for $\mu_B$ so that
$d(\mun_{A_N}, \mu_A) +d(\mun_{B_N}, \mu_B) = o(\sqrt N^{-1}).$\\
Thus, combining Theorem \ref{Fourier}.2
and Theorem \ref{sh} imply
$$\int_0^{2\theta} R_{\mu_B\boxplus \mu_A}(v)dv=
\int_0^{2\theta} R_{\mu_A}(v)dv
+\int_0^{2\theta} R_{\mu_B}(v)dv.$$
Differentiating with respect to $\theta$
gives Corollary \ref{additive}. \hfill\xx

Since the $R$-transform is analytic in a neighbourhood
of the origin, this entails the famous 
additivity property
of the $R$-transform.
So, Theorem \ref{Fourier}
provides a new proof of this
property, independent 
of cumulant techniques. \\

As announced in the introduction, 
the first step will be to use a result 
of concentration for orthogonal matrices.\\

\subsection{Concentration of measure for orthogonal matrices}

In this section, we prove the first point of Theorem \ref{Fourier}
that relies on
the following lemma, 
which is a direct consequence of a theorem due to Gromov \cite{Gro}

\begin{lem}
\label{aux}[Gromov, \cite{Gro}, p. 128]
Let $M_N^{(1)}$ denote the Haar measure 
on the special orthogonal group
$SO(N)$. 
There exists a positive constant $c>0$, independent
of $N$, 
such that for any function $F :SO(N)\ra \RR$
 so that there is a real $||F||_\La$ such that, for any $U,U'\in SO(N)$
$$|F(U)-F(U')|\ppq ||F||_\La \left(\sum_{i,j=1}^N |u_{ij}-u'_{ij}|^2
\right)^{\frac 1 2}, $$
then, for any $\e>0$,
$$M_N^{(1)}\left( \left|F(U)-\int F(U)dM_N^{(1)}(U)\right|\ge \e\right)
\le e^{-c N ||F||_\La^{-2} \e^2}.$$
\end{lem}
\textbf{ Proof of lemma \ref{aux} : }\\
 In \cite{Gro}, the author  prove such a 
lemma using the
fact that the Ricci curvature of $SO(N)$
is of order\footnote{In \cite{Gro}, it is reported that
the Ricci curvature is given by $N/4$ whereas 
J.C Sikorav and Y. Ollivier reported to
us that it is in fact $(N-2)/2$.} $N$ , and their 
result
 holds when  $F$ is Lipschitz with respect to the standard
bivariant metric which measures the length of the
geodesic  in $SO(N)$ between two elements
$U,U'\in SO(N)$. This distance   is of course
greater than the length of the geodesic in 
the whole space of matrices, given by the Euclidean 
distance, so that Lemma \ref{aux} is a direct consequence of
\cite{Gro}.
\hfill\xx

To prove Theorem
 \ref{Fourier}.1, %Ader
 we now apply our result 
with $F$ given by $F(U_N)=\frac 1 N \log I_N(\theta, A_N+U_NBU_N^*)$.
To get (\ref{Conc}), we have to check that this $F$ satisfies
the hypotheses of Lemma \ref{aux}. 
i.e. that $F$ is Lipschitz.

We have, for any matrices $W$, $\tilde W$
in $M_N := \{W \in \MNC; \, WW^* \ppq 1\}$,
\begin{multline*}
\left| \frac 1 N \log I_N(\theta, A_N+WB_NW^*)-\frac 1 N \log I_N(\theta, A_N+\tilde WB_N\tilde W^*)\right|\\
\ppq  2\theta ||B||_\infty 
\sup_{||v||= 1} \langle v, |W-\tilde W| v\rangle
\ppq 2\theta ||B||_\infty 
\left(\sum_{i,j=1}^N |w_{ij}-\tilde w_{ij}|^2\right)^{1\over 2}.
\end{multline*}
Moreover, if $T$ is for example the transformation changing the first column vector
$U_1$ of the matrix $U$ into $-U_1$,
$O(N)=SO(N)\sqcup T(SO(N)).$
Note that
$$F(TU) %={1\over N}\log I_N(\theta, A_N+(TU_N)B_N(TU_N)^*)
={1\over N}\log I_N(\theta, T^*A_NT+U_NB_N(U_N)^*).$$
Now, 
if we set
$E_N=A_N+U_NBU_N^*$ and $E_N'=T^*A_NT+U_NBU_N^*$,
we easily see that
$$d(\mun_{E_N}, \mun_{E_N'})\le {1\over N}\tr |E_N'-E_N|\le {2
||A||_\infty\over
N}.$$
Hence, 
 Lemma \ref{cont}.3 implies that
$$\d_N=\sup_{U\in SO(N)}
| F(U)- F(TU)|\ra 0 \mbox{ as }
N\ra\infty$$
Since
$$\int_{O(N)} F(U)dm_N^{(1)}(U)={1\over 2}\int_{SO(N)} F(U)dM_N^{(1)}(U)
+{1\over 2}\int_{SO(N)} F(TU)dM_N^{(1)}(U),$$
 we deduce that
$$\left|
\int_{O(N)} F(U)dm_N^1(U)-\int_{SO(N)} F(U)dM_N^1(U)\right|\le \d_N
.$$
Thus,  Lemma \ref{aux} implies that for $\e>0$
\begin{equation}\label{blork}
M_N^{(1)}\left( \left|F(U)-\int_{O(N)} 
F(U)dm_N^{(1)}(U)\right|\ge \e +\d_N\right)
\le e^{-c N ||F||_\La^{-2} \e^2}
\end{equation}
and similarly for $F(TU)$ so that
$$m_N^{(1)}\left( \left|F(U)-\int_{O(N)} F(U)dm_N^{(1)}(U)\right|\ge \e +\d_N\right)
\le e^{-c N ||F||_\La^{-2} \e^2}, $$
what gives Theorem \ref{Fourier}.1. \hfill\xx

\subsection{Exchanging  integration with the logarithm}
We are now seeking to establish the second point of Theorem \ref{Fourier}.
By Jensen's inequality,
$$\E[\log  I_N(\theta, A_N+V_N B_N (V_N)^*)]
\ppq \log \E[ I_N(\theta, A_N+V_N B_N (V_N)^*)]$$
so that we only need here  to
prove the converse inequality.

The whole idea to get it
 is contained
in the following

\begin{lem}
\label{control}
For any uniformly bounded sequence
of matrices $(A_N,B_N)_{N\in\N}$
and $\theta$ small enough, 
there exists 
a finite constant $C(A,B,\theta)$ such
that  for $N$ large enough
$${\E[ I_N(\theta, A_N+V_N B_N (V_N)^*)^2]\over 
\E[  I_N(\theta, A_N+V_N B_N (V_N)^*)]^2}\ppq C(\theta, A,B)$$
\end{lem}
Let us conclude the \textbf{proof of 
Theorem \ref{Fourier}.2} before proving this lemma.\\
 Hereafter, $\e>0$ is fixed. 
We introduce the event
$$\Aa =\left\{  I_N(\theta, A_N+V_N B_N (V_N)^*) \pgq {1\over 2}
\E[ I_N(\theta, A_N+V_N B_N (V_N)^*)]\right\}$$
Following \cite{Tal}, we have, 
if $I_N := I_N(\theta, A_N+V_N B_N(V_N)^*)$ that
$$ \EE[I_N] = \EE[I_N \mathbf 1_{\Aa^c}] + \EE[I_N \mathbf 1_{\Aa}]
\ppq \frac 1 2 \EE[I_N] + \EE[I_N^2]^{\frac 1 2 } \P(\Aa)^{\frac 1 2 }  $$
so that 
\begin{center} $ \displaystyle \frac{1}{4 C(A,B,\theta)} \ppq \P(\Aa).$\\
\end{center}
Furthermore, let $$t=
{1\over N}  \log \E\left[{1\over 2} I_N(\theta, A_N+V_N B_N
(V_N)^*)\right]-{1\over N}\E[\log I_N(\theta, A_N+V_N B_N(V_N)^*)]$$
We can assume that $t\pgq \d_N$ ($\d_N$ being given
in (\ref{blork})) since otherwise we are done.
We then  get  by (\ref{blork}) that  for any $t\ge\d_N$ and 
$N$ large enough,
$$\P(\Aa)\ppq \P\Bigl(  \frac 1 N \log I_N(\theta, A_N+UB_NU^*)
-m_N^{(1)}\Bigl( \frac 1 N \log I_N(\theta,
A_N+UB_NU^*)
\Bigr)\pgq t\Bigr) 
\ppq e^{-cN(t-\d_N)^2}$$
 with $c'=c(2|\theta|||B||_\infty)^{-2}$.
As a consequence,
$$\frac{1}{4C(A,B,\theta)}\ppq e^{-c'N(t-\d_N)^2}, 
\quad  \textrm{ so that }\quad 
t \ppq  \d_N +\left({1\over c' N}\log (4C(A,B,\theta))\right)^{1\over
2}
.$$
Hence, since $\d_N$ goes to zero with $N$,
$$
\lim_{N\ra\infty}\left({1\over N} \log \E\left[{1\over 2} I_N(\theta, A_N+V_N B_N
(V_N)^*)\right]  
-{1\over N}\E[\log I_N(\theta, A_N+V_N B_N(V_N)^*)]\right)
=0$$
which completes the proof of Theorem \ref{Fourier}.2.
 \hfill \xx

We go back to the \textbf{proof of Lemma \ref{control}.}
Observe first that
\begin{multline*}
L_N(\theta,A,B):=\E[ I_N(\theta, A_N+V_N B_N (V_N)^*)^2]\\
= \int e^{\theta N \left( (UAU^*)_{11}+(\tilde U A\tilde U^*)_{11}
+(UV_N B(UV_N)^*)_{11}+
(\tilde UV_N B(\tilde U V_N)^*)_{11}\right)}
dm_N^{(1)\,\otimes 3}(U,\tilde U,V_N)\\
= \int e^{\theta N \left( (UAU^*)_{11}+(\tilde U A\tilde
U^*)_{11}+ (VBV^*)_{11} +(\tilde U U^* V B V^*
U\tilde U^*)_{11}\right)  }
dm_N^{(1)\,\otimes 3}(V,U,\tilde U) 
\end{multline*}
where 
we used that $m_N^{(1)}$ is
invariant by the action of the orthogonal group.
We shall now prove that 
$L_N(\theta,A,B)$ factorizes.
The proof requires sharp estimates of 
spherical integrals. We already got the kind
of estimates we need in section \ref{tcl}. The ideas here will be 
very similar although the calculations will be more involved.\\

To rewrite $L_N(\theta,A,B)$ in a more proper way, the key observation 
is that, if we consider the column vector  $W:=(V^*
U\tilde U^*)_{1}$ then 
$\langle V_1,W \rangle= \langle U_1,\tilde U_1 \rangle$ so that we have the decomposition
$$W= \langle U_1,\tilde U_1 \rangle V_1+(1-|\langle U_1,\tilde U_1 \rangle|^2)^{{1\over 2}}V_2$$
with $(V_1,V_2)$ orthogonal and distributed 
uniformly on the sphere.\\
Therefore,
\begin{eqnarray*}
L_N(\theta,A,B)& = &
\E\left\lbk
\exp\{N\theta(F^N_1+F^N_2+F^N_3+F^N_4+F^N_5)\}\right\rbk
\end{eqnarray*}
with
\begin{eqnarray*}
F^N_1&=& \langle U,AU \rangle \\
F^N_2&=& \langle \tilde U,A\tilde U\rangle \\
F^N_3&=&(1+\langle U,\tilde U\rangle ^2)\langle V_1, B V_1\rangle \\
F^N_4&=& 2(1-|\langle U,\tilde U\rangle |^2)^{1\over 2}
\langle U,\tilde U\rangle  \langle V_1,BV_2\rangle \\
F^N_5&=& (1-\langle U,\tilde U\rangle ^2)\langle V_2,BV_2\rangle \\
\end{eqnarray*}
where $U$, $\tilde U$ are two independent vectors following 
the uniform law on the sphere
of radius $\sqrt{N}$ in $\RR^N$ and 
$V_1$, $V_2$ are the two
first column vectors 
of a matrix $V$ following $m_N^{(1)}$, $U$, $\tilde U$
and $V$ being 
independent.\\
\indent We now adopt the same strategy as in section \ref{tcl} 
to show that the $F_i$'s will become asymptotically independent
(or negligible). More precisely, we use again Fact
\ref{gaussien} and recall that we can write 
$\displaystyle U = \frac{g^{(1)}}{\|g^{(1)}\|}$, 
$\displaystyle \tilde U = \frac{g^{(2)}}{\|g^{(2)}\|}$,
$\displaystyle V_1 = \frac{g^{(3)}}{\|g^{(3)}\|}$
and $\displaystyle V_2 = \frac{G}{\|G\|}$ with
$\displaystyle G = g^{(4)} - \frac{\langle g^{(3)}, g^{(4)}\rangle }{\|g^{(4)}\|^2} g^{(3)}$ where
$g^{(1)}$, $g^{(2)}$, $g^{(3)}$ and $g^{(4)}$ are 4 i.i.d standard Gaussian vectors. We now set
for $i=1,2,3,4$, with $\l_j^{(i)}$ the eigenvalues 
of $A$ for $i=1$ or $2$ and of $B$ for $i=3$ or $4$,
$v_i= R_{\mu_A}(2\theta)$
for $i=1$ or $2$, $v_i= R_{\mu_B}(2\theta)$
for $i=3$ or $4$,
$$\hat U_i^N= {1\over N}\sum_{j=1}^N (g^{(i)}_j)^2-1,
\textrm{ and }
\hat V^N_i={1\over N}\sum_{j=1}^N \l_j^{(i)} (g^{(i)}_j)^2-v_i$$
Moreover, we let for $i=1$ or $2$,
$$\hat W^N_i={1\over N} \sum_{j=1}^N \l_j^{(i)} g^{(2i-1)}_j g^{(2i)}_j
\textrm{ and }
\hat Z^N_i= {1\over N} \sum_{j=1}^N  g^{(2i-1)}_j g^{(2i)}_j.$$

Under the Gaussian measure, all these quantities are going to zero almost surely
and we can localize $L_N$ as we made it in section \ref{short},
that is to say restrict the integration to the event
$\displaystyle A^\prime_N := \left\{ \hat U_i^N, \, \hat V^N_i, \, 
\hat W^N_i, \, \hat Z^N_i \textrm{ are } o( N^{-\frac 1 2 +\k})
\right\}$, for any $\k >0$. We then express  the $F_i$'s as function of these
variables and on $A^\prime_N$ we expand them till $o(N^{-1})$.
For example, on $A^\prime_N$,
$$ F_1={\hat V_1^N+v_1\over \hat U^N_1+1}
=  v_1+(\hat  V_1^N-v_1  \hat U^N_1)
-\hat U^N_1(\hat  V_1^N-v_1\hat U^N_1)+o(N^{-1})$$
and all the calculations go the same way so that
we get that the full
second order in $\sum_i F_i$ 
is
$$
\Xi^N=-\sum_{i=1}^4\hat U^N_i(\hat  V_i^N-v_i\hat U^N_i)
+2(\hat Z^N_1-\hat Z^N_2)\hat W^N_2- 2v_2\hat Z^N_2 \hat Z^N_1
+ 2v_2 (\hat Z^N_2)^2 $$

Now, as before, we consider the shifted probability measure $P_N$ 
(which contains all the first order term 
above) under which $(\tilde g^{(i)})_{i=1,\ldots, 4}$ defined by
$\tilde g^{(i)}_j = \sqrt{1+ 2\t v_i - 2 \t \l_j^{(i)}} g^{(i)}_j$ are i.i.d
standard Gaussian vectors.\\
Under $P_N$, the $(\hat U^N_i,\hat V^N_i)_{1\ppq i\ppq 4}$
are still independent with the same
law than for the
one dimensional
case.
 Moreover,
we see
that for $i=1,2,3,4$, $j=1,2$,
$$\lim_{N\ra\infty}
N\E[\hat U^N_i \hat Z^N_j]=0,\quad \lim_{N\ra\infty}
N\E[\hat U^N_i \hat W^N_j]=0
 .$$
Similarly, $(\hat Z^N_i,\hat W^N_i)_{i=1,2}$
are asymptotically uncorrelated.
Moreover, with $\mu_1=\mu_A$ and
$\mu_2=\mu_B$,
\begin{eqnarray*}
\lim_{N\ra\infty}
N\E[\hat W^N_i \hat Z^N_i]&=&\int {x\over (1+2\theta(v_i-x))^2} d\mu_i(x)
\\
\lim_{N\ra\infty}
N\E[(\hat W^N_i )^2]&=&\int {x^2\over ( 1+2\theta(v_i-x))^2} d\mu_i(x)
\\
\lim_{N\ra\infty}
N\E[(\hat Z^N_i )^2]&=&\int {1\over ( 1+2\theta(v_i-x))^2} d\mu_i(x).
\\
\end{eqnarray*}
Thus, with $G_i^N= \t v_i-{1\over  2N}\sum_{j=1}^N\log
(1-2\theta\l_j^{(i)}+2 \t v_i)$ and if the Gaussian integral 
is well defined, we have
\begin{multline*}
L_N(\theta, A,B)={e^{2N G_1^N+2NG_2^N}\over 
 \det(K_A) \det(K_B)}\\
\int \exp\{2\theta(\hat z_1-\hat z_2)\hat w_2- 2 v_2\theta\hat z_2\hat z_1
+ 2v_2\theta (\hat z_2)^2\}\prod_{i=1,2}dP_i(\hat w_i,\hat z_i)(1+o(1))
\end{multline*}
with  $P_i$  the law of two Gaussian 
variables with covariance matrix
$${R_i \over 2}=\left( \begin{array}{ll}
\int {1\over ( 1+2\theta(v_i-x))^2} d\mu_i(x)&\int 
{x\over ( 1+2\theta(v_i-x))^2} d\mu_i(x)\\
\int 
{x\over ( 1+2\theta(v_i-x))^2} d\mu_i(x)&\int 
{x^2\over ( 1+2\theta(v_i-x))^2} d\mu_i(x)\\
\end{array}\right)$$
and $K_A$ and $K_B$ as defined in (\ref{defK}) if we replace 
$\mu_E$ therein respectively by $\mu_A$ or $\mu_B$.\\
We now integrate on the variables $(\hat z_2,\hat w_2)$
so that the Gaussian computation gives
$$L_N(\theta, A,B)={e^{2N G_1^N+2NG_2^N}\over 
 \det(K_A) \det(K_B)^{3\over 2}}
\int \exp\{
\theta^2 \langle e,K_B^{-1} e \rangle \hat z_1^2 \}
dP_1(\hat z_1,\hat w_1)(1+o(1))$$
with $e=(-v_2, 1)$.
To show that the remaining integral is finite
it is enough to check that
$$ -2 \t^2 \langle e,K_B^{-1} e \rangle + \mbox{var} \hat z_1 \pgq 0,$$
at least for $\t$ small enough. But we can check that
$\theta^2 \langle e,K_B^{-1} e \rangle \approx \theta^2 \sigma_2$,
with $\s_2=\int x^2d\mu_B(x)$ whereas the variance 
of $\hat z_1$ is of order $1$.\\
This  finishes to
prove that for sufficiently small $\theta$'s 
there exists a finite constant $C(\theta,A,B)$
such that
 $$L_N(\theta, A,B)={e^{2N G_1^N+2NG_2^N}\over 
 \det(K_A) \det(K_B)}
C(\theta,A,B)(1+o(1))$$
Since on the other hand
we
have seen in section \ref{tcl} that
$$I_N(\theta,A)={e^{NG_1^N}\over \det{K_A}^{1\over 2}}(1+o(1))
\textrm{ and }
I_N(\theta,B)={e^{NG_2^N}\over \det{K_B}^{1\over 2}}(1+o(1)),$$
we have proved Lemma \ref{control}.\hfill \xx

\medskip

\section{Appendix}\label{app}
In this Appendix,
we clarify the derivation of the
central limit theorem of Theorems \ref{precise1}
and \ref{rg1complexe} and  Lemma \ref{control}.
We follow 
 the ideas
of \cite{bolt}, where only sums of i.i.d entries
$N^{-1}\sum_{i=1}^N  x_i$  were
considered rather than ponderated sums
$N^{-1}\sum_{i=1}^N \l_i x_i$.
 We consider the case of Theorem 
 \ref{rg1complexe} which is the most complicated;
\begin{equation}
\label{f11}
I_N(\t, E_N)
= \prod_{i=1}^N\sqrt{\z_i} \, \, e^{N\theta v}
\int \exp\left\{ N\theta{\g_N(v\g_N-\hat \g_N)\over
1+\g_N}\right\} e^{ -{1\over 2}\sum_{i=1}^N
g_i^2}\prod_{i=1}^N dg_i,
\end{equation}
where we recall that $\z_i := (1+2\t v - 2 \t \l_i)$,
$\g_N := \frac 1 N \sum_{i=1}^N \z_i g_i^2 - 1$ and $\hat \g_N = 
 \frac 1 N \sum_{i=1}^N \l_i \z_i g_i^2 - v$.
We denote
$$J_N(\t, E_N)=\sqrt{2\pi}^{-N}\int \exp\left\{ 
N\theta{\g_N(v\g_N-\hat \g_N)\over
1+\g_N}\right\} e^{ -{1\over 2}\sum_{i=1}^N
g_i^2}\prod_{i=1}^N dg_i.$$
The idea is the following :

\begin{itemize}
\item The first step is to derive a large deviation 
principle for $(\g_N, \hat \g_N)$ under 
the Gibbs measure
$$\mu_N^\theta(dg)=J_N(\t, E_N)^{-1}\exp\left\{ N\theta{\g_N(v\g_N-\hat \g_N)\over
1+\g_N}\right\} \prod_{i=1}^N P(dg_i).$$
As we showed that the unique minimizer is zero, it
entitles us to write
$$J_N(\t, E_N)=(1+\d(\e,\e',N))J_N^{\e,\e'}(\t, E_N)$$
with
$$J_N^{\e,\e'}(\t, E_N)=
\int_{|\g_N|\le \e, |\hat \g_N|\le\e'}
\exp\left\{ N\theta{\g_N(v\g_N-\hat \g_N)\over
1+\g_N}\right\} \prod_{i=1}^N P(dg_i)$$
where $\d(\e,\e',N)$ goes to zero 
as $N$ goes to infinity for any $\e,\e'>0$.

\item Let us assume that we can take above 
$\e=M/\sqrt{N},\e'= M'/\sqrt{N}$ with
$\d(M\sqrt{N}^{-1},M'\sqrt{N}^{-1} ,N)$
going to zero as $N$  and then $M,M'$ go to infinity. On the set
$\{|\g_N|\le N^{-\frac{1}{2}}M, |\hat \g_N|\le N^{-\frac{1}{2}}M'\}$,
$$f( \sqrt{N} \g_N,\sqrt{N} \hat\g_N)
= N\theta{ \g_N(v \g_N- \hat \g_N)\over
1+\g_N}= N\theta{\g_N(v \g_N- 
\hat \g_N)} +
O((M+M')^3 N^{-\frac{1}{2}})$$
and $f( \sqrt{N}\g_N,
\sqrt{N} \hat\g_N)$ is uniformly bounded. Further, the law
of $( N^{\frac{1}{2}} \g_N,
N^{\frac{1}{2}} \hat\g_N)$ converges under $P^{\otimes N}$ towards 
a two-dimensionnal complex Gaussian process with covariance matrix $
K'(\theta)$.
Hence, we can apply dominated convergence theorem
to see that 
$$\lim_{N\ra\infty}
\int_{|\g_N|\le N^{-\frac{1}{2}}M, 
|\hat \g_N|\le N^{-\frac{1}{2}}M'}
\exp\left\{ N\theta{\g_N(v\g_N-\hat \g_N)\over
1+\g_N}\right\} \prod_{i=1}^N P(dg_i)$$
$$
=(2\pi)^{-2}\mbox{det}(K'(\theta))^{-\frac{1}{2}}\int_{|x|\le M,|y|\le M'}
e^{ \theta x(vx-y)-\frac{1}{2} <(x,y), K'(\theta)^{-1} (x,y)> }
dxdy.$$
In the proof of Theorem \ref{precise1}, we established that
the bilinear
form $x,y\ra
\theta x(vx-y)-\frac{1}{2} <(x,y), K'(\theta)^{-1} (x,y)>$
 is strictly negative for $|\theta|$ small enough, 
therefore we can now let $M,M'$ going to infinity
to obtain a limit.

\item To see that we can take $\e=M/\sqrt{N},\e'=M'/\sqrt{N} $,
we can simplify the argument by recalling that
the spherical integral does not depend on $\gamma_N$.
Therefore,
$$
(1-P^{\ot N}(\e\ge |\gamma_N|\ge M\sqrt{N}^{-1}))
J^{\e,\e'}_N(\t, E_N)
=J^{MN^{-\frac12},\e'}_N(\t, E_N)$$
But, $\sqrt{N}\gamma_N= G_N^1+iG_N^2$ has, under $P^{\ot N}$,
sub-Gaussian exponential moments since
$$\E[e^{a G_N^j}]=\prod_{i=1}^N [ (1-2a\sqrt{N}^{-1}
\zeta_j(\l_i))^{-\frac{1}{2}}e^{- a\sqrt{N}^{-1}\zeta_j(\l_i)}]
\le e^{c a^2}$$
for some finite
constant $c$ which only depends
on a  uniform bound on the $\zeta_j(\l_i)$,
where we recall that $\zeta_j(\l_i) = \Re \zeta_i$ if $j=1$ and 
$\zeta_j(\l_i) = \Im \zeta_i$ if $j=2$. 
By Chebychev's inequality, we therefore conclude
that for $M$ big enough,
$$P^{\ot N}(|\gamma_N|\ge M\sqrt{N}^{-1})\le e^{-\frac{c}{8} M^2}.$$
Finally let us consider
$$J_N^{M,M',\e'}= \int_{|\g_N|\le  M\sqrt{N}^{-1},
 M'\sqrt{N}^{-1}\le |\hat \g_N|\le\e'}
\exp\left\{ N\theta{\g_N(v\g_N-\hat \g_N)\over
1+\g_N}\right\}\prod_{i=1}^N P(dg_i).$$
Clearly,  we find a finite constant $C$ (depending on $\theta$
and $\e'$) such that 
$$|J_N^{M,M',\e'}|\le e^{CM^2} \int_{|\g_N|\le  M\sqrt{N}^{-1},
 M'\sqrt{N}^{-1}\le |\hat \g_N|\le\e'}
\exp\left\{ C M|\sqrt{N}\hat \g_N|\right\} dP^{\ot N}(g).$$
Again, $\sqrt{N}\hat \g_N$ has sub-Gaussian tail
so that we find $C'>0$ so that 
$$|J_N^{M,M',\e'}|\le e^{(C+\frac{C^2}{C'}) M^2-C'(M')^2}.$$
Now, by
the previous point,
we know that
\begin{eqnarray*}
I(\theta,\mu_E)&=&\lim_{M,M'\ra\infty}
\lim_{N\ra\infty}
\int_{|\g_N|\le N^{-\frac{1}{2}}M, 
|\hat \g_N|\le N^{-\frac{1}{2}}M'}
e^{\left\{ N\theta{\g_N(v\g_N-\hat \g_N)\over
1+\g_N}\right\}} \prod_{i=1}^N P(dg_i)\\
&=& 
(2\pi)^{-2}\mbox{det}(K'(\theta))^{-\frac{1}{2}}
\int
e^{ \theta x(vx-y)-\frac{1}{2} <(x,y), K'(\theta)^{-1} (x,y)> }
dxdy
 \\
\end{eqnarray*}
exists and moreover goes to one
as $\theta$ goes to zero.
Hence, for $|\theta|$
small enough,
this term dominates $J_N^{M,M',\e'}$ for $N,M,M'$ large
enough ( $M'\gg M$) and we conclude
that
$$\lim_{N\ra\infty} J_N(\t, E_N)=\lim_{M,M'\ra\infty}
\lim_{N\ra\infty} J_N^{MN^{-\frac{1}{2}},M'N^{-\frac{1}{2}}}=
I(\theta,\mu_E).$$
Of course, this strategy only
requires non-degeneracy of the minimum
and $I(\theta,\mu_E)\neq 0$. 
In the setting of Theorem \ref{precise1},
this is verified on the whole interval
 $2\theta\in \HE([\l_{\rm min},\l_{\rm max}]^c)$.
In Lemma \ref{control}, we can also apply it by
noting that $L_N(\theta,A,B)$ does
not depend on \linebreak 
$(||g^{(1)}||, ||g^{(2)}||, ||g^{(3)}||,
||G||)$ to localize these quantities and proceed.

\end{itemize}

\nn
\nn
{\bf Acknowledgments :}
We are very grateful to O.~Zeitouni 
for helpful discussion at the beginning of this
work, which in particular allowed us to
obtain the second order correction in the full
high temperature region. We would like also to
thank P.~\'Sniady for many useful comments 
during this research. We thank Y.~Ollivier
 for pointing out \cite{Gro}
and showing us how Lemma \ref{aux} could
be deduced, which simplified a lot the argument.
Finally, we are also very grateful to an anonymous referee
whose careful reading and useful comments helped
us to improve the coherence of the paper.

\bibliography{rgfarxiv}

\begin{thebibliography}{10}

\bibitem{BA}
{\sc Ben~Arous, G.}
\newblock Methods de {L}aplace et de la phase stationnaire sur l'espace de
  {W}iener.
\newblock {\em Stochastics 25}, 3 (1988), 125--153.

\bibitem{BDG}
{\sc Ben~Arous, G., Dembo, A., and Guionnet, A.}
\newblock Aging of spherical spin glasses.
\newblock {\em Probab. Theory Related Fields 120}, 1 (2001), 1--67.

\bibitem{BGR}
{\sc Bercu, B., Gamboa, F., and Rouault, A.}
\newblock Large deviations for quadratic forms of stationary {G}aussian
  processes.
\newblock {\em Stochastic Process. Appl. 71}, 1 (1997), 75--90.

\bibitem{bolt}
{\sc Bolthausen, E.}
\newblock Laplace approximations for sums of independent random vectors.
\newblock {\em Probab. Theory Relat. Fields 72}, 2 (1986), 305--318.

\bibitem{Bor}
{\sc Borel, E.}
\newblock Sur les principes de la th\'eorie cin\'etique des gaz.
\newblock {\em Annales de l'\'Ecole Normale Sup\'erieure 23\/} (1906), 9--32.

\bibitem{Co}
{\sc Collins, B.}
\newblock Moments and cumulants of polynomial random variables on unitary
  groups, the {I}tzykson-{Z}uber integral, and free probability.
\newblock {\em Int. Math. Res. Not.}, 17 (2003), 953--982.

\bibitem{BP}
{\sc Collins, B., and Sniady, P.}
\newblock Integration with respect to the haar measure on unitary, orthogonal
  and symplectic group.
\newblock {\em Preprint, http://arxiv.org/abs/math-ph/0402073\/} (2004).

\bibitem{ADN}
{\sc D'Aristotile, A., Diaconis, P., and Newman, C.~M.}
\newblock Brownian motion and the classical groups.
\newblock In {\em Probability, statistics and their applications: papers in
  honor of Rabi Bhattacharya}, vol.~41 of {\em IMS Lecture Notes Monogr. Ser.}
  Inst. Math. Statist., Beachwood, OH, 2003, pp.~97--116.

\bibitem{DZ}
{\sc Dembo, A., and Zeitouni, O.}
\newblock {\em Large deviations techniques and applications}, second~ed.,
  vol.~38 of {\em Applications of Mathematics}.
\newblock Springer-Verlag, New York, 1998.

\bibitem{Gro}
{\sc Gromov, M., and Milman, V.~D.}
\newblock A topological application of the isoperimetric inequality.
\newblock {\em Amer. J. Math. 105}, 4 (1983), 843--854.

\bibitem{G}
{\sc Guionnet, A.}
\newblock First order asymptotics of matrix integrals ; a rigorous approach
  towards the understanding of matrix models.
\newblock {\em Comm. Math. Phys.\/} (2003).

\bibitem{GM1}
{\sc Guionnet, A., and Ma{\"{\i}}da, M.}
\newblock Character expansion method for the first order asymptotics of a
  matrix integral.
\newblock To appear in \textsl{Probab. Theory Related Fields}, 2005.

\bibitem{GZ}
{\sc Guionnet, A., and Zeitouni, O.}
\newblock Large deviations asymptotics for spherical integrals.
\newblock {\em J. Funct. Anal. 188}, 2 (2002), 461--515.

\bibitem{HC1}
{\sc Harish-Chandra}.
\newblock Differential operators on a semisimple {L}ie algebra.
\newblock {\em Amer. J. Math. 79\/} (1957), 87--120.

\bibitem{HC2}
{\sc Harish-Chandra}.
\newblock Fourier transforms on a semisimple {L}ie algebra. {I}.
\newblock {\em Amer. J. Math. 79\/} (1957), 193--257.

\bibitem{ItZ}
{\sc Itzykson, C., and Zuber, J.~B.}
\newblock The planar approximation. {II}.
\newblock {\em J. Math. Phys. 21}, 3 (1980), 411--421.

\bibitem{Ji}
{\sc Jiang, T.}
\newblock How many entries of a typical orthogonal matrix can be approximated
  by independent normals ?
\newblock Preprint, 2003.

\bibitem{karlin}
{\sc Karlin, S., and McGregor, J.}
\newblock Coincidence probabilities.
\newblock {\em Pacific J. Math. 9\/} (1959), 1141--1164.

\bibitem{K2}
{\sc Kazakov, V.}
\newblock Solvable matrix models.
\newblock In {\em Random matrix models and their applications}, vol.~40 of {\em
  Math. Sci. Res. Inst. Publ.} Cambridge Univ. Press, Cambridge, 2001,
  pp.~271--283.

\bibitem{MPR}
{\sc Marinari, E., Parisi, G., and Ritort, F.}
\newblock Replica field theory for determinstic models. {II}. {A} non-random
  spin glass with glassy behaviour.
\newblock {\em J. Phys. A 27}, 23 (1994), 7647--7668.

\bibitem{Ma}
{\sc Matytsin, A.}
\newblock On the large-{$N$} limit of the {I}tzykson-{Z}uber integral.
\newblock {\em Nuclear Phys. B 411}, 2-3 (1994), 805--820.

\bibitem{Sag}
{\sc Sagan, B.~E.}
\newblock {\em The symmetric group}.
\newblock The Wadsworth \& Brooks/Cole Mathematics Series. Wadsworth \&
  Brooks/Cole Advanced Books \& Software, Pacific Grove, CA, 1991.
\newblock Representations, combinatorial algorithms, and symmetric functions.

\bibitem{Tal}
{\sc Talagrand, M.}
\newblock Concentration of measure and isoperimetric inequalities in product
  spaces.
\newblock {\em Inst. Hautes \'Etudes Sci. Publ. Math.}, 81 (1995), 73--205.

\bibitem{VoiR}
{\sc Voiculescu, D.}
\newblock Addition of certain noncommuting random variables.
\newblock {\em J. Funct. Anal. 66}, 3 (1986), 323--346.

\bibitem{Vo1}
{\sc Voiculescu, D.}
\newblock The analogues of entropy and of {F}isher's information measure in
  free probability theory. {I}.
\newblock {\em Comm. Math. Phys. 155}, 1 (1993), 71--92.

\bibitem{Zel}
{\sc Zelditch, S.}
\newblock Macdonald's identities and the large {$N$} limit of {$YM\sb 2$} on
  the cylinder.
\newblock {\em Comm. Math. Phys. 245}, 3 (2004), 611--626.

\bibitem{ZZJ}
{\sc Zinn-Justin, P., and Zuber, J.-B.}
\newblock On some integrals over the {${\rm U}(N)$} unitary group and their
  large {$N$} limit.
\newblock {\em J. Phys. A 36}, 12 (2003), 3173--3193.
\newblock Random matrix theory.

\bibitem{ZV}
{\sc Zvonkin, A.}
\newblock Matrix integrals and map enumeration: an accessible introduction.
\newblock {\em Math. Comput. Modelling 26}, 8-10 (1997), 281--304.
\newblock Combinatorics and physics (Marseilles, 1995).

\end{thebibliography}
\bibliographystyle{acm}

\end{document}